\documentclass[12pt,leqno,fleqn,epsfig]{article}
\usepackage{amssymb, epsfig, amsmath, amsthm}
\usepackage{mathrsfs}       

\newcommand{\R}{\mathbb{R}}

\newcommand{\Z}{\mathbb{Z}}

\newcommand{\N}{\mathbb{N}}

\newcommand{\PP}{\mathbb{P}}

\newcommand{\mes}{\operatorname{\rm meas}}    
\newcommand{\esssup}{\operatorname*{ess\,sup}}

\newcommand{\supp}{\operatorname*{supp}}

\newcommand{\const}{\operatorname*{const}}

\newcommand{\be}{\begin{equation}}
\newcommand{\ee}{\end{equation}}
\newcommand{\bea}{\begin{eqnarray}}
\newcommand{\eea}{\end{eqnarray}}
\newcommand{\bean}{\begin{eqnarray*}}
\newcommand{\eean}{\end{eqnarray*}}

\newcommand{\var}{\varepsilon}

\newcommand{\intl}{\int\limits}

\newcommand{\Beweisende}{\rule{0.2cm}{0.2cm}}

\newcounter{secnum}

\newtheorem{thm}{Theorem}[section]

\newtheorem{lem}[thm]{Lemma}

\theoremstyle{definition}

\newtheorem{defin}[thm]{Definition}
\newtheorem{rem}[thm]{Remark}

\title{Existence of discretely self-similar solutions to the Navier-Stokes equations  for initial value in $ L^2_{ loc}(\R^{3})$} 
 
\author{Dongho Chae$^*$  and J\"{o}rg Wolf $^\dagger$\\
\ \\
 $*$Department of Mathematics\\
Chung-Ang University\\
 Seoul 156-756, Republic of Korea\\
 e-mail: dchae@cau.ac.kr\\
and \\
$\dagger$Department of Mathematics\\
Humboldt University Berlin\\
Unter den Linden 6, 10099 Berlin, Germany\\
e-mail: jwolf@math.hu-berlin.de}
\date{}
\begin{document}

\maketitle
\begin{abstract}
We prove the existence of a forward discretely self-similar  solutions to the Navier-Stokes equations in  $ \R^{3}\times (0,+\infty)$ for a discretely  self-similar initial velocity belonging to $ L^2_{ loc}(\R^{3})$.   
  \noindent{\bf
keywords:} Navier-Stokes  equations, existence, discretely self-similar solutions 

\end{abstract}

\section{Introduction}
\label{sec:-1}
\setcounter{secnum}{\value{section} \setcounter{equation}{0}
\renewcommand{\theequation}{\mbox{\arabic{secnum}.\arabic{equation}}}}

In this paper we study the existence of forward discretely self-similar (DSS) solutions to the Navier-Stokes equations in 
$Q= \R^{3}\times (0, +\infty)$
\begin{align}
\nabla \cdot  u &=0,  
\label{1.1}
\\
\partial _t u  + (u \cdot \nabla ) u -\Delta u  &= - \nabla \pi,
\label{1.2}
\end{align}
with the initial condition 
\begin{equation}
u  = u_0 \quad  \text{on}\quad  \R^{3}\times \{0\}.
\label{1.3}
\end{equation}
Here $ u(x,t)= (u_1(x,t), u_2(x,t), u_3(x,t))$ denotes the velocity of the fluid, and $ u_0(x)= (u_{ 0,1}(x), u_{ 0,2}(x), u_{ 0,3}(x))$, 
while $ \pi $ stands for the pressure.  In case $ u_0 \in L^2(\R^{3})$ with $ \nabla \cdot u_0=0$ in the sense of distributions the global in time existence of weak solutions to \eqref{1.1}--\eqref{1.3} , which satisfy the global energy inequality for almost all $ t\in (0,+\infty)$
\begin{equation}
 \frac{1}{2} \| u(t)\|^2_{ 2} +  \intl_{0}^{t} \| \nabla u(s)\|^2_2  ds \le  \frac{1}{2} \| u_0\|^2_2
\label{1.4}
\end{equation}
has been proved  by Leray \cite{Leray1934}. On the other hand, the important questions of regularity and 
uniqueness of solutions to \eqref{1.1}--\eqref{1.3} are still open. The first significant  results in this direction have been established by Scheffer 
\cite{Scheffer1976} and later  by Caffarelli, Kohn, Nirenberg   \cite{CaKoNi1982} for solutions $ (u,\pi )$ that also  satisfy the 
following local energy inequality for almost all $ t\in (0,+\infty)$ and for all nonnegative $ \phi \in 
C^{\infty}_{\rm c}(Q)$
\begin{align}
& \frac{1}{2}  \intl_{\R^{3}} | u( t)|^2 \phi(x,t) dx +  \intl_{0}^{t}   \intl_{\R^{3}} |\nabla  u|^2 \phi dx  ds
 \cr
&\qquad \le   \frac{1}{2} \intl_{0}^{t}   \intl_{\R^{3}}  | u|^2  \Big(\frac{\partial}{\partial t} +\Delta\Big)  \phi dx  ds
+  \frac{1}{2} \intl_{0}^{t}   \intl_{\R^{3}}  (| u|^2 + 2\pi) u\cdot \nabla \phi   dx  ds.
\label{1.5}
\end{align}
On the other hand, the space $ L^2(\R^{3})$  excludes homogenous spaces of degree $ -1$ belonging  to the scaling invariant class. In fact we observe that 
$ u_\lambda(x,t) = \lambda u(\lambda x, \lambda ^2t)$ solves the Navier-Stokes equations with initial velocity 
$ u_{ 0, \lambda }(x)= \lambda u_0(\lambda x)$, for any $ \lambda >0$.  This suggests to study of the Navier-Stokes system 
for initial velocities in a homogenous space $ X$ of degree $ -1$, which means that  $ \| v\|_{ X}= \| v_\lambda \|_X$ 
for all $ v\in X$.  Koch and Tataru proved  in \cite{KoTa2001} that $ X= BMO^{ -1}$  is the largest possible space with scaling invariant norm which guarantees  well-posedness under smallness condition. On the contrary, for self-similar (SS) initial data fulfilling  $ u_{ 0,\lambda }=u$ for all $ \lambda >0$  a natural space seems to be  $ X= L^{3,\infty}(\R^{3})$. This space is  embedded into the space 
$ L^2_{ uloc}(\R^{3})$,  which contains  uniformly local square integrable functions.  Obviously, possible solutions to the Navier-Stokes equations with $ u_0 \in L^2_{ uloc}(\R^{3})$ 
do not satisfy the global energy equality, rather the local energy inequality in the sense of Caffarelli-Kohn-Nirenberg.  Such solutions 
are called local Leray solutions. The existence of global in time local Leray solutions  has been proved by  Lemari\`e-Rieusset in 
 \cite{LemRie2002} (see also in \cite{KiSe2007} for more details). This concept has been used by Bradshaw and Tsai \cite{BradTsai2015}  for the construction of a discretely self-similar ($ \lambda $-DSS, $ \lambda >1$) local Leray solution  for a $ \lambda $-DSS initial velocity $ u_0 \in L^{ 3,\infty}(\R^{3})$. This result generalizes the previous results  of Jia and \v{S}ver\'ak \cite{JiaSve2014} 
 concerning  the existence of SS  local Leray solution, and the result  by Tsai in \cite{Tsai2014}, which proves   the existence of a $ \lambda $-DSS Leray solution for $ \lambda $ near $ 1$. 
 However, for the $ \lambda $-DSS initial data it would be more natural to assume $ u_0 \in L^2_{ uloc}(\R^{3})$ instead  $ L^{ 3, \infty}(\R^{3})$.  In general, such initial value does not belong to 
$ L^2_{ uloc}(\R^{3})$ and therefore it does not belong to the Morrey class $ M^{ 2,1}$, rather to the weighted space $ L^{ 2}_{ k}(\R^{3})$ of all $ v\in L^2_{ loc}(\R^{3})$ such that $ \frac{v}{(1+| x|^k)}\in L^2(\R^{3})$ for all 
$ \frac{1}{2} < k <+\infty$.

\hspace{0.5cm}
 Since the authors in   \cite{BradTsai2015}  work on the existence of  periodic solutions to the 
time dependent Leray equation a certain spatial decay is necessary which can be ensured for initial data in $ L^{ 3,\infty}(\R^{3})$. 
On the other hand, applying the local $ L^2$ theory it would 
be more natural to assume $ u_0 \in L^2(B_\lambda  \setminus B_1)$ only. 
As explained  in  \cite{BradTsai2015} their method even breaks down for initial data in the Morrey class $ M^{ 2,1}(\R^{3})$, which is a much smaller subspace of  $ L^2_{ loc}(\R^{3})$.   By using an entirely different 
method we are able to construct a global weak solutions for such DSS initial data.  
  
  \hspace{0.5cm}
In the present paper we introduce a new notion of a local Leray solution   satisfying a local energy inequality with projected pressure.  To the end, we provide  the  notations of function spaces which will be use in the sequel. 
By $ L^s(G ), 1 \le s \le \infty$ we denote the usual Lebesgue spaces. The usual Sobolev spaces are denoted by $ W^{k,\, s}(G)$ 
and $ W^{k,\, s}_0(G), 1 \le s \le +\infty, k\in \N$. The dual of $ W^{k,\, s}_0(G)$ will be denoted by $ W^{-k,\, s'}(G)$, where $ s'= \frac{s}{s-1},  1 < s < +\infty$.  For a general space of vector fields $ X$  the subspace of solenoidal fields  will be denoted by $ X_{ \sigma }$. In particular, the space of solenoidal smooth fields with compact support is denoted by $ C^{\infty}_{\rm c, \sigma }(\R^{3})$. In addition we define energy space 
\[
 V^2 (G\times (0,T))= L^\infty(0,T; L^2(G))\cap L^2(0,T; W^{1,\, 2}(G)),\quad  0<T \le +\infty.  
\]

\hspace{0.5cm}
We now recall the definition of the local pressure projection $ E^{ \ast}_{G}: W^{-1,\, s}(G) 
\rightarrow  W^{-1,\, s}(G)$ for a given bounded $ C^2$-domain $ G \subset \R^{3}$,  introduced in 
\cite{Wolf2016} based on the unique solvability of the steady Stokes system (cf. \cite{GaSiSo1994}).  
More precisely, for any  $ F\in W^{-1,\, s}(G)$ there exists a unique pair $ (v, p)\in W^{1,\, s}_{ 0, \sigma }(G)\times L^s_0(G)$ 
 which solves weakly the steady Stokes system 
 \begin{equation}
 \begin{cases}
 \nabla \cdot v =0\quad  \text{ in}\quad  G,\quad  -\Delta v + \nabla p = F\quad  \text{ in}\quad  G,
  \\[0.3cm]
 v=0   \text{ on}\quad  \partial G.
 \end{cases}
 \label{1.5b}
 \end{equation}
 
 Here $ W^{1,\, s}_{ 0, \sigma }(G)$ stands for closure of $ C^{\infty}_{\rm c, \sigma }(\R^{3})$  with respect to the norm in  $ W^{1,\, s}(G)$, while $ L^s_0(G)$ denotes the subspace of $ L^s(G)$ with vanishing average. 
Then we set $E^{ \ast}_G(F):=  \nabla p $, where $ \nabla p $ denotes the gradient function in $ W^{-1,\, s}(G)$  defined as 
\[
 \langle \nabla p, \varphi \rangle = -\intl_{G} p \nabla \cdot \varphi dx, \quad  \varphi \in W^{1,\, s'}_0(G). 
\]

\begin{rem}
\label{rem1.2}
From the existence and uniqueness of weak solutions $ (v,p)$ to \eqref{1.5b} for given for any  $ F\in W^{-1,\, s}(G)$ it follows that 
\begin{equation}
\| \nabla v\|_{ s, G} + \| p\|_{ s, G} \le c \| F\|_{ -1,s, G}. 
\label{1.6a}
\end{equation} 
where $ c= \const $ depending on $ s$ and the geometric properties of $ G$,  and  depends only on $ s$  if $ G$ equals a ball or an annulus which due to the scaling properties of the Stokes equation.  
In case $ F$ is given by $ \nabla \cdot f$ for $f \in L^s(\R^{3})^9$ then \eqref{1.6a} gives 
\begin{equation}
 \| p\|_{ s, G} \le c \| f\|_{ s, G}. 
\label{1.6b}
\end{equation}  
According to the estimate $ \| \nabla p\|_{ -1,s, G}  \le \| p\|_{ s, G}$, and using \eqref{1.6b}, we see that the operator $ E^{ \ast}_G$ 
is bounded in $ W^{-1,\, s}(G)$. Furthermore, as $ E^{ \ast}_G(\nabla p)=\nabla p$ for all $ p\in L^s_{ 0}(G)$ we see that $ E^{ \ast}_G$ defines a projection.

\hspace{0.5cm}
2. In case $ F\in L^s(G)$, using the canonical embedding  $ L^s(G) \hookrightarrow  W^{-1,\, s}(G)$, by the aid of  elliptic regularity 
 we get $ E^{ \ast}_G(F)= \nabla p \in L^s(G)$ together with the estimate 
\begin{equation}
\| \nabla p\|_{s, G} \le c \| F\|_{ s, G}, 
\label{1.6c}
\end{equation}
where the  constant in \eqref{1.6c}  
depends only on $ s$ and $ G$. In case  $ G$  equals a ball or an annulus this constant depends only on $ s$ 
 (cf. \cite{GaSiSo1994} for more details). Accordingly the restriction of  $ E^{ \ast}_G$ to the Lebesgue space $ L^s(G)$ appears to be  
a projection in $L^s(G)$. This projection will be denoted still by $ E^{ \ast}_G$.

\end{rem}

 \begin{defin}[Local Leray solution with projected pressure] 
  \label{def1.1}
 Let $ u_0\in L^2_{loc}(\R^{3})$. A vector function
 $ u\in L^2_{loc, \sigma  }(\R^{3}\times [0, +\infty))$ is called a {\it local Leray solution to \eqref{1.1}--\eqref{1.3} with projected pressure}, if 
 for any bounded $ C^2$  domain $ G \subset \R^{3}$ and $ 0<T<+\infty$
 \begin{itemize}
  \item[1.]   $ u\in V^2_\sigma (G\times (0,T))\cap C_{ w}([0,T]; L^2(G))$.
  \item[2.]  $u $ is a distributional solution to \eqref{1.2} , i.\,e. for every $ \varphi \in C^{\infty}_{\rm c}(Q)$ with $ \nabla \cdot \varphi =0$  
  \begin{equation}
  \intl  \hspace*{-0.2cm}\intl  _{ \hspace*{-0.5cm} Q}- u\cdot \frac{\partial \varphi }{\partial t} - u \otimes u : \nabla \varphi + \nabla u : \nabla \varphi dxdt =0. 
 \label{1.6}
 \end{equation}
 \item[3.]  $  u(t) \rightarrow u_0$ in $ L^2(G)$ as $ t \rightarrow 0^+$.
 \item[4.] The following local energy inequality with  projected pressure  holds  for  every nonnegative $ \phi \in C^{\infty}_{\rm c}(G\times (0,+\infty))$, and for almost every $ t\in (0,+\infty)$ 
 \begin{align}
 &   \frac{1}{2} \intl_{G} | v_G(t)|^2 \phi dx +    \intl_{0}^{t}\intl_{G} | \nabla v_{G}|^2  \phi dx    ds
\cr
 & \quad  \le  \frac{1}{2}  \intl_{0}^{t} \intl_{G} | v_{G}|^2 \Big(\Delta + \frac{\partial }{\partial t}\Big) \phi  + 
 | v_G|^2 u\cdot \nabla \phi )  dx    ds
 \cr
 & \qquad +  
 \intl _{0}^{t}\intl_{G} (u \otimes v_{G}) :\nabla ^2 p_{h, G} \phi dx dt   +\intl_{0}^{t} \intl_{G} p_{ 1,G } v_{G} \cdot \nabla \phi dxds
 \cr
 &\qquad \qquad \qquad +  
 \intl_{0}^{t}\intl_{G} p_{2,G } v_{G}\cdot \nabla \phi dxds,
 \label{1.7}
\end{align} 
 where $ v_{G}= u + \nabla p_{ h, G }$, and 
\begin{align*}
\nabla p_{ h,G } &= -E^{ \ast}_{G} (u),
\\
\nabla p_{ 1,G } &= -E^{ \ast}_{ G} ((u\cdot \nabla ) u), \quad \nabla p_{2,G } = E^{ \ast}_{G} (\Delta u). 
\end{align*} 

  \end{itemize} 

 \end{defin}
 
 \begin{rem}
 1. Note that due to $ \nabla \cdot u=0$ the pressure $ p_{ h, G}$ is harmonic, and thus smooth in $ x$. 
 Furthermore, as it has been proved in \cite{Wolf2016} the pressure gradient  $ \nabla p_{ h,G}$ is continuous 
 in $ G\times [0,+\infty)$. 
 
 \hspace{0.5cm}
 2. The notion of local suitable weak solution to the Navier-Stokes equations satisfying the local energy inequality 
 \eqref{1.7}  has been  introduced in \cite{Wolf2015c}.  As it has been shown there such solutions enjoy the same partial regularity properties as the usual suitable weak solutions in the  Caffarelli-Kohn-Nirenberg theorem. 
 
 \end{rem}

\hspace{0.5cm}
Our main result is the following 

\begin{thm}
\label{thm1.4} For any $ \lambda $-DSS initial data $u_0\in L^2_{loc, \sigma }(\R^{3}) $ there exists at least 
one local Leray  solution with projected pressure  $ u\in L^2_{ loc, \sigma }(\R^{3}\times [0, +\infty))$ to the Navier-Stokes equations 
\eqref{1.1}--\eqref{1.3} in the sense of Definition\,\ref{def1.1}, which is discretely self-similar.  

\end{thm}

 \section{Solutions of the linearized problem with initial velocity in $ L^2_{ \lambda - DSS}$}
 \label{sec:-4}
 \setcounter{secnum}{\value{section} \setcounter{equation}{0}
 \renewcommand{\theequation}{\mbox{\arabic{secnum}.\arabic{equation}}}}
 
 Let $ 1< \lambda < +\infty$ be fixed.  For $ f: \R^{3} \rightarrow \R^{3}$ we denote 
  $ f_\lambda(x):= \lambda f(\lambda x), x\in \R^{3}$. For a time dependent function $ f: Q \rightarrow \R^{3}$  we denote $ f_\lambda(x,t):= \lambda f(\lambda x, \lambda ^2 t), (x, t)\in \R^{3}\times (0,+\infty)$.  We now define for $ 1 \le s \le +\infty$
 \begin{align*}
 L^s_{ \lambda - DSS}(\R^{3}) &:= \Big\{u\in L^1_{ loc}(\R^{3}) \,\Big|\, u \in L^s(B_\lambda  \setminus B_1), u_\lambda =u \,\,\text{ a.\,e. in $ \R^{3}$}\Big\}, 
 \\
 L^s_{ \lambda - DSS}(Q) &:= \Big\{u\in L^1_{ loc}(Q) \,\Big|\, u \in L^s(Q_\lambda  \setminus Q_1), u_\lambda =u\,\,\text{ a.\,e. in $ Q$}\Big\}. 
 \end{align*}
 Here $ B_r$ stands usual ball in $ \R^{3}$ with center $ 0$ and radius $ r>0$, while $ Q_r= B_r \times (0, r^2)$. 
 
 \hspace{0.5cm}
 In the present section we consider the following linearized problem in $ Q$
 \begin{align}
\nabla \cdot  u &=0,  
\label{4.1}
\\
\partial _t u  + (b \cdot \nabla ) u -\Delta u  &= - \nabla \pi
\label{4.2}
\end{align}
with the initial condition 
\begin{equation}
u  = u_0 \quad  \text{on}\quad  \R^{3}\times \{0\},
\label{4.3}
\end{equation}
where $ u_0$ belongs to $ L^2_{ \lambda -DSS}(\R^{3})$  with $ \nabla\cdot u_0=0$,   and $ b\in L^s_{ \lambda -DSS}(Q)$, $ 3 \le s \le 5$ with $ \nabla\cdot b=0$ both  in the sense of distributions. 
We give the following notion of a local solution  with projected pressure for the linear system \eqref{4.1}, \eqref{4.2}.  
 
 \begin{defin}[Local solution with projected pressure to the linearized problem] 
  \label{def4.1}
 Let $ u_0\in L^2_{loc, \sigma }(\R^{3})$ and let $ b\in L^3_{loc, \sigma}(\R^{3}\times [0,+\infty))$. A vector function
 $ u\in L^2_{loc, \sigma }(\R^{3}\times [0, +\infty))$ is called a {\it local  solution to \eqref{4.1}--\eqref{4.3} with projected pressure}, if 
 for any bounded $ C^2$  domain $ G \subset \R^{3}$ and $ 0<T<+\infty$ the following conditions are satisfied
 \begin{itemize}
  \item[1.]   $ u\in V^2(G\times (0,T))\cap C_{ w}([0,T]; L^2(G))$.
  \item[2.]  $u $ is a distributional solution to \eqref{4.2} , i.\,e. for every $ \varphi \in C^{\infty}_{\rm c}(Q)$  with $ \nabla \cdot \varphi =0$
  \begin{equation}
  \intl  \hspace*{-0.2cm}\intl  _{ \hspace*{-0.5cm} Q}- u\cdot \frac{\partial \varphi }{\partial t} - b \otimes u : \nabla \varphi + \nabla u : \nabla \varphi dxdt =0. 
 \label{4.0a}
 \end{equation}
 \item[3.]  $  u(t) \rightarrow u_0$ in $ L^2(G)$ as $ t \rightarrow 0^+$.
 \item[4.] the following local energy inequality with projected pressure  holds  for  every nonnegative $ \phi \in C^{\infty}_{\rm c}(G\times (0,+\infty))$, and for almost every $ t\in (0,+\infty)$ 
  
 \begin{align}
 &   \frac{1}{2} \intl_{G} | v_G(t)|^2 \phi dx +    \intl_{0}^{t}\intl_{G} | \nabla v_{G}|^2  \phi dx    ds
\cr
 & \quad  \le  \frac{1}{2}  \intl_{0}^{t} \intl_{G} | v_{G}|^2 \Big(\Delta + \frac{\partial }{\partial t}\Big) \phi  + 
 | v_G|^2 b\cdot \nabla \phi )  dx    ds
 \cr
 & \qquad +  
 \intl _{0}^{t}\intl_{G} (b \otimes v_{G}) :\nabla ^2 p_{h, G} \phi dx dt   +\intl_{0}^{t} \intl_{G} p_{ 1,G } v_{G} \cdot \nabla \phi dxds
 \cr
 &\qquad \qquad \qquad +  
 \intl_{0}^{t}\intl_{G} p_{2,G } v_{G}\cdot \nabla \phi dxds
 \label{4.0b}
\end{align} 
 where $ v_{G}= u + \nabla p_{ h, G }$, and 
\begin{align*}
\nabla p_{ h,G } &= -E^{ \ast}_{G} (u),
\\
\nabla p_{ 1,G } &= -E^{ \ast}_{ G} ((b\cdot \nabla ) u), \quad \nabla p_{2,G } = E^{ \ast}_{G} (\Delta u). 
\end{align*} 

  \end{itemize} 

 \end{defin}

 \begin{thm}
 \label{thm4.1}
Let  $ b\in L^3_{ \lambda -DSS}(Q)\cap L^{ \frac{18}{5}}(0,T; L^3(B_1))$, $ 0< T < +\infty$,  
with $ \nabla \cdot b=0$ in the sense 
of distributions. Suppose that $ b\in L^3_{ loc}(0, \infty; L ^\infty(\R^{3}))$. 
  For every $ u_0\in L^2_{ \lambda -DSS}(\R^{3})$ with $ \nabla \cdot u_0=0$ in the sense of distributions,   there exists a unique local 
   solution  with projected pressure $ u\in L^2_{loc, \sigma }(\R^{3}\times [0, +\infty))$  to \eqref{4.1}--\eqref{4.3}  according to 
 Definition\,\ref{def4.1}
such that for any $ 0< \rho < +\infty$ and $ 0< T <+\infty$ it holds 
 \begin{align}
 & u\in L^3_{\lambda -DSS}(Q),
 \label{4.1d}
 \\
 & u \in  C([0,T]; L^2(B_\rho )),
 \label{4.1a}
 \\
 & \| u\|_{ L^\infty(0, T; L^2(B_{\rho ^{  \frac{1}{4}}}))} + \| \nabla u\|_{ L^2(B_{ \rho^{\frac{3}{5} } }\times (0,T))}
\le C_0 K_0\Big(\rho ^{ \frac{1}{2}}+ |\!|\!|b|\!|\!|^3 \max\{ T^{ \frac{13}{18}}, T^{  \frac{1}{2} }\}\Big),
 \label{4.2b}
 \\
 &\| u\|_{ L^4(0,T; L^3(B_1))} \le C_0 K_0\Big(1+ |\!|\!|b|\!|\!|^3 \max\{ T^{ \frac{13}{18}}, T^{  \frac{1}{2} }\}\Big),
 \label{4.2bb}
  \end{align}
  where $ K_0 := \| u_0\|_{ L^2(B_1)}$ and  $ |\!|\!|b|\!|\!|= 
  \| b\|_{ L^{ \frac{18}{5}}(0,T; L^3(B_1))}$, 
  while $ C_0>0$ denotes a  constant depending on $ \lambda $ only.

 \end{thm}
 
 \hspace{0.5cm}
 Before turning to the proof of Theorem\,\ref{4.1}, we  show the existence and uniqueness  of 
  weak solutions to the  linear system \eqref{4.1}-- \eqref{4.3} for $ L^2_\sigma $ initial data.

 \begin{lem}
 \label{lem4.2}
  Let  $ b\in L^3_{ \lambda -DSS}(Q)\cap L^{ \frac{18}{5}}(0,T; L^3(B_1))$, $ 0< T < +\infty$ with $ \nabla \cdot b=0$ in the sense 
of distributions. Suppose that $ b\in L^3_{ loc}(0, \infty; L ^\infty(B_1))$. 
  For every $ u_0\in L^2_{\sigma }(\R^{3})$ there exists a unique 
 weak solution  $ u\in V^2_\sigma (Q)\cap  C([0,+\infty); L^2(\R^{3}))$  to \eqref{4.1}--\eqref{4.3}, which satisfies the global 
 energy equality for all $ t\in [0, +\infty)$
 \begin{equation}
  \frac{1}{2} \| u(t)\|^2_{ 2} +  \intl_{0}^{t} \intl_{ \R^{3}} | \nabla u|^2  dxds =  \frac{1}{2} \| u_0\|^2_{ 2}.
 \label{4.30}
 \end{equation}  
 
 \end{lem}   
 
{\bf  Proof}: {\it 1. Existence}: By using standard linear theory of parabolic systems we easily get the existence of a weak solution 
$ u\in V^2(Q)\cap C_w([0,+\infty); L^2(\R^{3}))$ to \eqref{4.1}--\eqref{4.3} which satisfies the global energy inequality  
for almost all $ t\in (0, +\infty)$ 
 \begin{equation}
  \frac{1}{2} \| u(t)\|^2_{ 2} +  \intl_{0}^{t} \intl_{ \R^{3}} | \nabla u|^2  dxds \le   \frac{1}{2} \| u_0\|^2_{ 2}.
 \label{4.31}
 \end{equation}  
It is well known that such solutions have the property 
\begin{equation}
 u(t) \rightarrow u_0  \quad  \text{{\it in}}\quad  L^2(\R^{3})\quad  \text{{\it as}}\quad  t \rightarrow 0^+.  
\label{4.32}
\end{equation}
On the other hand, from the assumption of the Lemma it follows that  for all $ t_0\in (0, T)$
\[
\| b u\|_{ L^2( \R^{3}\times (t_0, T))} \le \| b\|_{ L^{ \infty}(\R^{3}\times (t_0, T))}\| u_0\|_{2}.  
\]
Accordingly, $ u\in C((0, T]; L^2(\R^{3}))$ and for all $ t_0\in (0,T]$ and $ t\in [t_0, T]$ the following energy equality holds true 
\begin{equation}
  \frac{1}{2} \| u(t)\|^2_{ 2} +  \intl_{t_0}^{t} \intl_{ \R^{3}} | \nabla u|^2  dxds =  \frac{1}{2} \| u(t_0)\|^2_{ 2}. 
 \label{4.33}
 \end{equation}  
  Now letting $ t_0 \rightarrow 0$ in \eqref{4.33}, and observing \eqref{4.32}, we are led to \eqref{4.30}.     
 
 \hspace{0.5cm}
 By a similar argument, making use of \eqref{4.32} we easily prove the local energy inequality \eqref{4.0b}.   
 
 \hspace{0.5cm}
 {\it 2. Uniqueness}: Let $ v \in V^2_\sigma (Q)$ be a second solution to \eqref{4.1}--\eqref{4.3} satisfying the global energy equality. 
 As we have seen above this solution belongs to $ C([0,+\infty); L^2(\R^{3}))$. Setting  $ w = u-v$, 
 by our  assumption on $ b$  it follows that   $ b \otimes w\in L^2(\R^{3}\times (t_0, T])$ for any  $ t_0 \in (0, T]$. Accordingly,  as above we get the following energy equality 
 \begin{equation}
  \frac{1}{2} \| w(t)\|^2_{ 2} +  \intl_{t_0}^{t} \intl_{ \R^{3}} | \nabla w|^2  dxds =  \frac{1}{2} \| w(t_0)\|^2_{ 2}. 
 \label{4.34}
 \end{equation}  
Verifying that $ w(t_0) \rightarrow 0$ in $ L^2(\R^{3})$ as $ t_0 \rightarrow 0^+$ from \eqref{4.34} letting  $ t_0 \rightarrow 0^+$ 
it follows that $ \| w(t)\|_{ 2}= 0 $ for all $ t\in [0,T]$. This completes the proof of the uniqueness.  \hfill \Beweisende

 \vspace{0.5cm}  
{\bf  Proof of Theorem\,\ref{thm4.1}}: Since $ u_0$ is $ \lambda $-DSS we have $ \lambda u_0(\lambda x) = u_0(x)$ for all $ x\in \R^{3}$. We define the extended annulus $ {\tilde A} _k = B_{ \lambda ^k}  \setminus B_{ \lambda ^{ k-3}}$, $ k\in \N$.  Clearly, $ B_1 \cup  (\cup _{ k=1}^\infty {\tilde A} _k) = \R^{3}$. There exists a partition of unity $ \{ \psi _k\}$ such that   $ \supp\psi _k\subset {\tilde A} _k$ for $ k\in \N$ and 
$\supp \psi _0 \subset  B_1 $, and  $ 0 \le \psi_k \le 1$,  $ | \nabla^2 \psi _k|+ | \nabla \psi _k|^2 \le \lambda ^{ -2k}$, $ k\in \N\cup \{ 0\}$. 
We set $ u_{ 0, k} = \PP (u_0 \psi _k )$, $ k\in \N \cup \{ 0\}$. Clearly, 
\begin{equation}
 u_0 = \sum_{k=0}^{\infty} u_{ 0,k},
 \label{4.2a}
\end{equation}
where the limit in \eqref{4.2a} is taken in the sense of  $ L^2_{ loc}(\R^{3})$. 

\hspace{0.5cm}
Let $ k\in \N\cup \{ 0\}$ be fixed.  Thanks to Lemma\,\ref{lem4.2} we get  a unique weak solution $ u_k \in V^2_{ \sigma }(Q)$ to the problem 
 \begin{align}
\nabla \cdot  u_k &=0 \quad  \text{in}\quad  Q,
\label{4.3a}
\\
\partial _t u_k  + (b \cdot \nabla ) u_k -\Delta u_k  &= - \nabla \pi_k  \quad  \text{in}\quad  Q,
\label{4.4}
\\
u_k &= u_{ 0,k} \quad  \text{on}\quad  \R^{3}\times \{0\},
\label{4.5}
\end{align}
satisfying the following global  energy equality for all $ t\in [0,+\infty)$
\begin{equation}
\frac{1}{2}\| u_{ k}(t)\|_2^2 +  \intl_{0}^{t} \intl_{  \R^{3}} | \nabla u_k| ^2 dx ds =\frac{1}{2}\| u_{0, k}\|_2^2. 
\label{4.6}
\end{equation}
 By using the transformation formula, we get 
 \begin{align}
\| u_{0, k}\|_2^2  & \le  \intl_{ \R^{3}} | u_0 \psi _k |^2 dx
\le  \intl_{{\tilde A} _k} | u_0  |^2 dx =  \lambda ^{3k} \intl_{{\tilde A} _1} 
| u_0(\lambda^k x)|^2 dx 
 \cr
 &= \lambda^k \intl_{{\tilde A} _1} | \lambda^k  u_0(\lambda^k x)|^2 dx
= \lambda^k \intl_{{\tilde A} _1} |  u_0(x)|^2 dx \le cK_0^2\lambda^{ k}. 
\label{4.6a}
 \end{align}
Combining \eqref{4.6} and \eqref{4.6a}, we are led to 
\begin{equation}
\| u_k\|^2_{ L^\infty(0,T; L^2)}+  \| \nabla u_k\|^2_{ L^2(0,T; L^2)}  \le  
cK_0^2\lambda^{ k}.
\label{4.6b}
\end{equation}  

\hspace{0.5cm}
Next, let $ \lambda ^{ \frac{3}{5}k } \le r < \rho \le \lambda ^{ \frac{3}{5}(k +1)} $ be arbitrarily chosen, but fixed. 
By introducing the local pressure we have 
\[
\frac{\partial v_{ k, \rho }}{\partial t} + (b\cdot \nabla ) u_k - \Delta v_{k, \rho } = - \nabla p_{1, k, \rho } - \nabla p_{ 2, k, \rho },
\]
 where $ v_{k, \rho }= u_k + \nabla p_{ h,k, \rho }$, and 
\begin{align*}
\nabla p_{ h,k, \rho } &= -E^{ \ast}_{B_{ \rho}} (u_k),
\\
\nabla p_{ 1,k, \rho } &= -E^{ \ast}_{ B_{ \rho}} ((b\cdot \nabla ) u_k), \quad \nabla p_{ 2,k, \rho } = E^{ \ast}_{ B_{ \rho} } (\Delta u_k). 
\end{align*} 
 The following local energy equality holds true for all $ \phi \in C^{\infty}_{\rm c}(B_{ \rho})$  and for all $ t\in [0,T]$, 
 \begin{align}
 &\frac{1}{2} \intl_{B_{ \rho} } | v_{k, \rho}(t)|^2 \phi^{ 6}  dx  +  \intl_{0}^{t} \intl_{B_{ \rho}} | \nabla v_{k, \rho}|^2  \phi^{6}dx    ds 
 \cr
 &= \frac{1}{2}  \intl_{0}^{t} \intl_{B_{ \rho} } | v_{k, \rho}|^2\Delta \phi^{6}  dx    ds + 
 \frac{1}{2}  \intl_{0}^{t} \intl_{B_{ \rho} } | v_{k, \rho}|^2 b\cdot \nabla \phi ^{ 6} dx    ds 
 \cr
 &\qquad\qquad  +
 \intl_{0}^{t} \intl_{B_{ \rho} } (b \otimes v_{k, \rho}) :\nabla ^2 p_{h, k, \rho} \phi ^{6}dx ds  
 +  \intl_{0}^{t} \intl_{B_{ \rho} } p_{ 1,k, \rho } v_{k, \rho} \cdot \nabla \phi^{ 6} dxds
 \cr
 & \qquad \qquad +  
 \intl_{0}^{t} \intl_{B_{ \rho} } p_{2,k, \rho } v_{k, \rho}\cdot \nabla \phi^{ 6} dxds+  \frac{1}{2}
  \intl_{B_{ \rho} }  | v_{0, k}|^2 \phi ^{ 6}  dx  
 \cr
 &= I + II+ III+ IV+ V+ VI.
  \label{4.7}
 \end{align}
 
 Let $ \phi\in  C^{\infty}_{\rm c}(\R^{3})$ denote a cut off function such that $ 0 \le \phi \le 1$ in $ \R^{3}$, 
 $ \phi \equiv 1$ on $ B_{r}$, $ \phi \equiv 0$ in $ \R^{3}  \setminus B_{\rho }$, and $ | \nabla ^2 \phi |+ | \nabla \phi |^2 
 \le c (\rho -r)^{ -  2}$ in $ \R^{3}$. 

\hspace{0.5cm}
Let $ m\in \N$ be chosen so that  $  \lambda ^{m-1} \le \rho < \lambda ^{m}$. Then we estimate 

 \begin{align*}
 \| b\|_{ L^3(B_\rho \times (0,T))}^3 &=  \lambda ^{ 5m} \intl_{0}^{T \lambda ^{ -2m}} \intl_{B_{ \rho \lambda ^{ -m}}}  | b(\lambda ^{ -m} x,  \lambda ^{ -2m} t)|^3  dxdt 
 \\
 &=\lambda ^{ 2m} \intl_{0}^{T \lambda ^{ -2m}} \intl_{B_{ \rho \lambda ^{ -m}}}  | b(x,  t)|^3  dxdt 
 \\
 & \le c\lambda ^{ 2m - \frac{1}{3}m} T^{ \frac{1}{6}}
 \| b\|_{ L^{ \frac{18}{5}}(0,T; L^3(B_1))}^3 \le c |\!|\!|b|\!|\!|^3\rho ^{ \frac{5}{3}}  T^{ \frac{1}{6}}, 
 \end{align*}
 where and hereafter the constants appearing in the estimates may depend on $ \lambda $.
 The above estimate together with 
 $ \rho^{ \frac{5}{3}} \le \lambda^{k+1 } $ yields
 \begin{equation}
  \| b\|_{ L^3(B_\rho \times (0,T))} \le   c |\!|\!|b|\!|\!|
  \lambda ^{ \frac{1}{3}k}  T^{ \frac{1}{18}}.  
 \label{4.7}
 \end{equation}
In what follows we extensively make use of the estimate for almost all $ t\in (0,T)$ 
\begin{equation}
\| \nabla p_{ h,k,\rho }(t)\|_{ L^2(B_\rho )}  \lesssim  \| u_k(t)\|_{ L^2(B_\rho )}, 
\label{4.7a}
\end{equation} 
which is an immediate consequence of \eqref{1.6c}.  In addition, we easily verify the inequality 
\begin{equation}
\| \nabla^2 p_{ h,k,\rho }(t)\|_{ L^2(B_\rho )}  \lesssim  \| \nabla u_k(t)\|_{ L^2(B_\rho )}.
\label{4.7b}
\end{equation} 
Indeed, observing that 
\[
 \nabla^2 p_{ h,k,\rho }(t) = \nabla (\nabla  p_{ h,k,\rho }(t)- u(t)_{ B_\rho })
= -\nabla E_{ B_\rho }^{ \ast} (u_k(t)- u_k(t)_{ B_\rho })
\]
by means of elliptic regularity along with the Poincar\'e inequality we get 
\begin{align*}
\| \nabla^2 p_{ h,k,\rho }(t) \|^2_{ L^2(B_\rho )} 
&\le c \rho ^{ -2}\| u_k(t)- u_k(t)_{ B_\rho }\|^2_{ L^2(B_\rho )} +  c \| \nabla u_k(t)\|^2_{ L^2(B_\rho )}
\\
&\le  c \| \nabla u_k(t)\|^2_{ L^2(B_\rho )}.
\end{align*}
Whence, \eqref{4.7b}.

\hspace{0.5cm}
(i) With the help of  \eqref{4.6b} we easily deduce that 
\begin{align*}
I  \le  c (\rho -r) ^{ - 2} \intl_{0}^{t} \intl_{ B_{ \rho} } |  u_{k}|^2   dxds
\le c K_0^2 (\rho -r) ^{ - 2}  \lambda ^{k}T. 
\end{align*}  

(ii) Next,  using H\"older's inequality and Young's inequality together with  \eqref{4.6b}, \eqref{4.7}, \eqref{4.7a}  and \eqref{4.7b}, we estimate      
 \begin{align*}
II &\le    (\rho -r) ^{ - 1}\intl_{0}^{t} \intl_{B_{ \rho} } | b| |  v_{k, \rho}|^2 \phi^{ 5}  dxds
  \\
  & \le c (\rho -r) ^{ - 1} T^{ \frac{1}{6}}
    \| b\|_{ L^3(B_\rho \times (0,T))} \| v_{ k,\rho } \phi ^{3}\|_{ L^{ \infty}(0,T;L^2)} 
  \| v_{ k,\rho }\phi ^{ 2}\|_{ L^2(0,T; L^6)}  
  \\
  & \le c (\rho -r) ^{ - 2} T^{ \frac{2}{3}}
    \| b\|_{ L^3(B_\rho \times (0,T))} \| v_{ k,\rho } \phi ^{3}\|_{ L^{ \infty}(0,T;L^2)} 
  \| u_k\|_{ L^\infty(0,T; L^2)}  
 \\
  &\qquad \qquad +c (\rho -r) ^{ - 1} T^{ \frac{1}{6}}
    \| b\|_{ L^3(B_\rho \times (0,T))} \| v_{ k,\rho } \phi ^{3}\|_{ L^{ \infty}(0,T;L^2)} 
  \| \nabla v_{ k,\rho }\phi ^{ 2}\|_{ L^2(0,T; L^2)}    
       \\
  & \le   c|\!|\!|b|\!|\!|K_0 (\rho -r) ^{ - 2} \lambda ^{ \frac{5}{6} k} T^{ \frac{13}{18} }    \| v_{k, \rho  }  \phi ^3 \|_{ L^\infty(0,T; L^2)}
  \\ 
  &\qquad +  c|\!|\!|b|\!|\!|(\rho -r) ^{ - 1} \lambda ^{ \frac{1}{3}k}  T^{ \frac{2}{9}}
   \| v_{k, \rho  } \phi ^{ 3}\|_{ L^\infty(0,T; L^2)} \| \nabla v_{ k, \rho } \phi ^{ 3}\|_{ L^2(0,T; L^2)}^{ \frac{2}{3}} \| \nabla v_{ k, \rho } \|_{ L^2(0,T; L^2(B_\rho ))}^{ \frac{1}{3}}
 \\
 & \le   c|\!|\!|b|\!|\!|K_0 (\rho -r) ^{ - 2} \lambda ^{ \frac{5}{6} k} T^{ \frac{13}{18} }    \| v_{k, \rho  }  \phi ^3 \|_{ L^\infty(0,T; L^2)}
  \\ 
  & \qquad  +  c|\!|\!|b|\!|\!| K_0^{ \frac{1}{3}}(\rho -r) ^{ - 1}  \lambda ^{ \frac{1}{2}k} T^{ \frac{2}{9}}
   \| v_{k, \rho  } \phi ^{ 3}\|_{ L^\infty(0,T; L^2)} \| \nabla v_{ k, \rho } \phi ^{ 3}\|_{ L^2(0,T; L^2)}^{ \frac{2}{3}} 
       \\
 &\le c  |\!|\!|b|\!|\!|^2K_0^2 (\rho -r) ^{ - 4} \lambda ^{ \frac{5}{3} k}T^{ \frac{13}{9}}  
  \\ 
  &\qquad   +  c|\!|\!|b|\!|\!|^6K_0^2(\rho -r) ^{ - 6} \lambda ^{ 3k} T^{ \frac{4}{3}}  + 
  \frac{1}{8}  \| v_{k, \rho  } \phi ^{ 3}\|^{2}_{ L^\infty(0,T; L^2)} 
  + \frac{1}{4} \| \nabla v_{ k, \rho } \phi ^{ 3}\|^2_{ L^2(0,T; L^2)}  
  \\ 
  & \le   cK_0^2 (\rho -r)^{ -3} \lambda ^{k} \max\{ T^{ \frac{13}{9}}, T\}
  + c|\!|\!|b|\!|\!|^6K_0^2(\rho -r) ^{ - 6} \lambda ^{ 3k} \max\{ T^{ \frac{13}{9}}, T\}
  \\
  &  \qquad \qquad + 
  \frac{1}{8}  \| v_{k, \rho  } \phi ^{ 3}\|^{2}_{ L^\infty(0,T; L^2)} 
  + \frac{1}{4} \| \nabla v_{ k, \rho } \phi ^{ 3}\|^2_{ L^2(0,T; L^2)}. 
 \end{align*}
 (iii) In what follows we make use the following estimates using the fact that $ p_{ h, k, \rho }$ is harmonic. 
 By using the identity 
  \[
\intl_{\R^{3}} | \nabla h|^2 \phi ^2 dx= \frac{1}{2} \intl_{\R^{3}} h^2 \Delta \phi ^2 dx 
\]
for any harmonic function $ h$ on $B_\rho $, and cut off function $ \phi \in C^{\infty}_{\rm c}(B_\rho )$, we get 
\begin{align}
\| \nabla ^3 p_{ h, k, \rho }(t) \phi ^3 \|_{ 2} \le c (\rho -r)^{ -1}\| \nabla ^2 p_{ h, k, \rho }(t) \phi ^2 \|_2 \le (\rho -r)^{ -2}\| \nabla  p_{ h, k, \rho }(t)  \|_{ 2, B_\rho }
\label{4.7d}
\end{align}
  By the aid of  Sobolev's inequality, together with \eqref{4.7d}, we get for almost every $ t\in (0,T)$ 
 \begin{align*}
 \| \nabla^2 p_{ h, k, \rho }(t)\phi^3\|_{6} &\le c (\rho - r)^{ -1} \| \nabla ^2 p_{ h, k, \rho }(t) \phi ^2\|_{ 2, B_\rho } + c \| \nabla ^3 p_{ h, k, \rho }(t) \phi ^3 \|_{ 2}
 \\
 &\le  c (\rho - r)^{ -1} \| \nabla ^2 p_{ h, k, \rho }(t)\phi ^2\|_{ 2, B_\rho }
 \\
 &\le  c (\rho - r)^{ -2} \| \nabla p_{h,  k, \rho  }(t)\|_{ 2, B_\rho }
 \\
 &\le  c (\rho - r)^{ -2} \|u_{k}(t)\|_{ 2, B_\rho }. 
 \end{align*}
 Integrating both sides of the above estimate, and  
 estimating the right-hand side of of the resultant inequality  by \eqref{4.6b}, we arrive at 
 \begin{equation}
 \| \nabla^2 p_{ h, k, \rho} \phi^3 \|_{L^2(0,T; L^6)} \le c (\rho - r)^{ -2} 
 T^{ \frac{1}{2}} K_0 \lambda ^{ \frac{1}{2}k}. 
 \label{4.7e}
 \end{equation}

 Arguing as above, and using \eqref{4.7e},   we find   
 \begin{align*}
 III & \le cT^{ \frac{1}{6}} \| b\|_{ L^3(0,T L^3(B_{ \rho } )} 
 \| v_k \phi ^{ 3}\|_{ L^\infty(0,T; L^2)} \| \nabla ^2 p_{h,k,\rho  } \phi ^{3}\|_{ L^2(0,T; L^6)}
 \\
 &\le cK_0 (\rho -r)^{ -2}T^{ \frac{2}{3} }   \lambda ^{ \frac{1}{2}k}\| b\|_{ L^3(0,T L^3(B_{\rho } )} 
 \| v_k \phi ^{ 3}\|_{ L^\infty(0,T; L^2)}
 \\
 &\le c|\!|\!|b|\!|\!|K_0 (\rho -r)^{ -2} \lambda ^{\frac{1}{2} k}T^{ \frac{13}{18} }    
 \| v_k \phi ^{ 3}\|_{ L^\infty(0,T; L^2)} 
 \\
 &\le c|\!|\!|b|\!|\!|^2K_0^2 (\rho -r)^{ -4} \lambda ^{k}T^{ \frac{13}{9}}+ \frac{1}{8}
 \| v_k \phi ^{ 3}\|^2_{ L^\infty(0,T; L^2)}.
   \end{align*}
 
 (iv) We now  going to estimate $ IV$. Using \eqref{1.6b}, and arguing similar as before,  we estimate   
 \begin{align*}
 IV &\le  c (\rho -r)^{ -1}\| p_{ 1, k, \rho }\|_{ L^{ \frac{6}{5}}(0,T; L^{ 2}(B_{ \rho} ))} \| v_{k, \rho} \phi ^{ 3}\|_{ L^6(0,t; L^2)}
 \\
 &\le c  (\rho -r)^{ -1} T^{ \frac{1}{6}} \| b u_k\|_{ L^{ \frac{6}{5}}(0,T; L^{ 2}(B_{ \rho} ))} 
  \|  v_{k, \rho}\phi ^{ 3}\|_{ L^\infty(0,T; L^2)}
 \\
 &\le c  (\rho -r)^{ -1} T^{ \frac{1}{6}}\| b \|_{L^3(0,T; L^{3}(B_\rho ))} \| u_k \|_{ L^2(0,T; L^{6}(B_\rho ))} 
  \|  v_{k, \rho}\phi ^{ 3}\|_{ L^\infty(0,T; L^2)}
  \\
 &\le c |\!|\!|b|\!|\!| (\rho -r)^{ -1} \lambda ^{ \frac{1}{3}k} T^{ \frac{2}{9}}
 \| u_k \|_{ L^2(0,T; L^{6}(B_\rho ))} 
  \|  v_{k, \rho}\phi ^{ 3}\|_{ L^\infty(0,T; L^2)}
  \\
 &\le c |\!|\!|b|\!|\!|K_0 (\rho -r)^{ -1}  \rho^{ -1} \lambda ^{ \frac{1}{3}k} T^{ \frac{13}{18}}   
    \|  u_{k}\|_{ L^\infty(0,T; L^2)}  \|  v_{k, \rho}\phi ^{ 3}\|_{ L^\infty(0,T; L^2)}
 \\
 & \qquad +
  c |\!|\!|b|\!|\!| (\rho -r)^{ -1}  \lambda ^{ \frac{1}{3}k} T^{ \frac{2}{9}} 
  \|  v_{k, \rho}\phi ^{ 3}\|_{ L^\infty(0,T; L^2)}
  \| \nabla u_k \|^{ \frac{1}{3}}_{ L^2(0,T; L^{2}(B_\rho ))}  \| \nabla u_k \|^{ \frac{2}{3}}_{ L^2(0,T; L^{2}(B_\rho ))}
 \\
 &\le c |\!|\!|b|\!|\!|K_0 (\rho -r)^{ -1}  \lambda ^{ \frac{7}{30}k} T^{ \frac{13}{18}}   
   \|  v_{k, \rho}\phi ^{ 3}\|_{ L^\infty(0,T; L^2)}
 \\
 & \qquad +
  c |\!|\!|b|\!|\!|K_0^{ \frac{1}{3}} (\rho -r)^{ -1}  \lambda ^{ \frac{1}{2}k} T^{ \frac{2}{9}} 
  \|  v_{k, \rho}\phi ^{ 3}\|_{ L^\infty(0,T; L^2)}\| \nabla u_k \|^{ \frac{2}{3}}_{ L^2(0,T; L^{2}(B_\rho ))} 
  \\
 &\le c |\!|\!|b|\!|\!|^2 K_0^2 (\rho -r)^{ -2}  \lambda ^{ \frac{7}{15}k} T^{ \frac{13}{9}} +c |\!|\!|b|\!|\!|^6K_0^2 (\rho -r)^{ -6}  \lambda ^{ 3k} T^{ \frac{4}{3}} 
  \\
 & \qquad  +   \frac{1}{8} 
 \|  v_{k, \rho}\phi ^{ 3}\|^2_{ L^\infty(0,T; L^2)} + \frac{1}{4}\| \nabla u_k \|^{2}_{ L^2(0,T; L^{2}(B_\rho ))}
 \\ 
  & \le (1+|\!|\!|b|\!|\!|^6)K_0^2(\rho -r) ^{ - 6} \lambda ^{ \frac{17}{5}k} \max\{ T^{ \frac{13}{9}}, T\}
 \\
 &\qquad \qquad  + 
  \frac{1}{8}  \| v_{k, \rho  } \phi ^{ 3}\|^{2}_{ L^\infty(0,T; L^2)} 
  + \frac{1}{4}\| \nabla u_k \|^{2}_{ L^2(0,T; L^{2}(B_\rho ))}.
     \end{align*}

(v)  Recalling the definition of $ p_{ 2, k, \rho }$, using \eqref{1.6b},   \eqref{4.6b} and Young's inequality, we get 
\begin{align*}
V &\le   c   (\rho -r)^{ -1} \| p_{ 2,k, \rho } \|_{L^2(0, T; L^2(B_{ \rho} )) } \|  v_{k, \rho} \phi^{ 3} \|_{ L^2(0,T; L^2)} 
\\
&\le c   (\rho -r)^{ -1} T^{ \frac{1}{2}} \bigg( \intl_{0}^{T} \intl_{B_{ \rho} } | \nabla u_k|^2  dx   dt \bigg)^{  \frac{1}{2} }
\| v_{k, \rho} \phi ^{ 3}\|_{ L^\infty(0,T; L^2)} 
\\
&\le c  K_0^2 (\rho -r)^{ -2} \lambda ^{ k} T + \frac{1}{8}\| v_{k, \rho} \phi ^{ 3}\|^2_{ L^\infty(0,T; L^2)}
\\
&\le c  K_0^2 (\rho -r)^{ -6} \lambda ^{ \frac{17}{5} k} T + \frac{1}{8}\| v_{k, \rho} \phi ^{ 3}\|^2_{ L^\infty(0,T; L^2)}.
\end{align*}

(vi) It only remains to evaluate  $ VI$.  Let $ k \ge 9$. Then $ \frac{3}{5} (k +1)\le  k-3$. Thus, $ \supp(\psi _k)\cap B_\rho = \emptyset$. In particular, $ \psi _k u_{0} =0$ in $ B_\rho $. This shows that,  almost  everywhere in $ B_\rho $ it holds 
\[
u_{ 0, k} = \PP( \psi _k u_{0}) - \psi _k u_{0}
\]
which is a gradient field. Accordingly, almost  everywhere in $ B_\rho $
\begin{align*}
 v_{ 0, k} &= u_{ 0, k } - E^{ \ast}_{ B_\rho }(u_{ 0, k }) =  u_{ 0, k} - u_{ 0,k}=0.  
\end{align*}
Hence 
\[
VI =0. 
\]
For $ k \le 8$    we find 
\[
VI \le \| u_{ 0,k}\|^2_{L^2(B_\rho )} \le c \sum_{k=0}^{8} \| u_0 \psi _k\|^2_{ L^2}
\le c\| u_{ 0}\|^2_{ L^2(B_{ \lambda ^8 })} \le c K_0^2. 
\]

We now insert the above estimates of $ I, \ldots, VI$ into the right-hand side of \eqref{4.7}. This gives 
\begin{align}
&\esssup_{t\in (0,T)} \intl_{B_\rho }  | v_{k, \rho}(t)|^2 \phi ^6 dx +  \intl_{0}^{T} \intl_{B_\rho }  | \nabla v_{k, \rho}|^2 \phi ^6 dxdt
\cr
&\qquad \le  cK_0 ^2 \max\{ 8-k, 0\} + c(1+|\!|\!|b|\!|\!|^6)K_0^2 \max\{ T^{ \frac{13}{9}}, T\}(\rho -r) ^{ - 6} \lambda ^{ \frac{17}{5}k}  
\cr
&\qquad \qquad + \frac{1}{4}\intl_{0}^{T} \intl_{B_\rho }  | \nabla u_{k}|^2dxdt. 
\label{4.8}
\end{align}
On the other hand, employing \eqref{4.7d} and \eqref{4.6b}  
\[
\intl_{B_\rho }  | \nabla^2 
p_{h, k, \rho}|^2 \phi ^6 dxdt 
\le c K_0^2 (\rho -r)^{ -2} \lambda^k T, 
\]
we estimate 
\begin{align}
&\intl_{0}^{T} \intl_{B_r }  | \nabla u_{k}|^2 dxdt
\cr
& \qquad \le  2\intl_{0}^{T} \intl_{B_\rho }  | \nabla v_{k, \rho}|^2 \phi ^6 dxdt + 2\intl_{0}^{T} \intl_{B_\rho }  | \nabla^2 
p_{h, k, \rho}|^2 \phi ^6 dxdt
\cr
&\qquad \le  2\intl_{0}^{T} \intl_{B_\rho }  | \nabla v_{k, \rho}|^2 \phi ^6 dxdt + 
c K_0^2 (\rho -r)^{ -2} \lambda ^{ k} T
\label{4.9}
\end{align}
 Combining \eqref{4.8} and \eqref{4.9}, we are led to 
 \begin{align}
  &\intl_{0}^{T} \intl_{B_r}  | \nabla u_{k}|^2 dxdt
\cr
&\qquad \le cK_0^2 \max\{ 8-k, 0\}+  c(1+|\!|\!|b|\!|\!|^6)K_0^2 \max\{ T^{ \frac{13}{9}}, T\}(\rho -r) ^{ - 6} \lambda ^{ \frac{17}{5}k}  
\cr
&\qquad \qquad + \frac{1}{2}\intl_{0}^{T} \intl_{B_\rho }  | \nabla u_{k}|^2dxdt. 
\label{4.10}
\end{align}
By virtue of a routine iteration argument from \eqref{4.10}  we get for all 
$ \rho \in [ \lambda ^{ \frac{3}{5}k}, 2\lambda ^{ \frac{3}{5}k}]$
 \begin{align}
 & \esssup_{t\in (0,T)} \intl_{B_{ \rho /2} }  | v_{k, \rho}(t)|^2  dx+  \intl_{0}^{T} \intl_{B_{ \rho/2 }}  | \nabla u_{k}|^2 dxdt
  \cr
 &\qquad \qquad \le cK_0^2 \max\{ 8-k, 0\}+  c(1+|\!|\!|b|\!|\!|^6)K_0^2 \max\{ T^{ \frac{13}{9}}, T\}\rho  ^{ - 6} \lambda ^{ \frac{17}{5}k}
 \cr
 &\qquad \qquad \le cK_0^2 \max\{ 8-k, 0\}+  c(1+|\!|\!|b|\!|\!|^6)K_0^2 \max\{ T^{ \frac{13}{9}}, T\}\lambda ^{ -\frac{1}{5}k}.
\label{4.11}
\end{align}
  
  \hspace{0.5cm}
  In addition, by using the mean value property of harmonic functions along with \eqref{4.6b},   we  estimate for almost all $ t\in (0,T)$ 
  \begin{align*}
\|\nabla p_{h, k, \rho}(t)\|^2_{ L^2(B_{ \lambda ^{  \frac{1}{4}k}})} 
&\le c \lambda ^{ \frac{3}{4} k} \| \nabla p_{ h,k, \rho}(t)\|^2_{ L^\infty(B_{ \rho /2} )}
  \\
&\le c 
\lambda ^{ -\frac{21}{20}k}\| \nabla p_{h, k, \rho}(t)\|^2_{ L^2(B_\rho )}
  \\
&\le c 
\lambda ^{ -\frac{21}{20}k}\| u_k\|^2_{ L^\infty(0,T; L^2(B_\rho) )}
\le c K_0^2\lambda ^{ - \frac{1}{20}k}.  
  \end{align*}
 \[
\]
Combining this estimate with \eqref{4.11}, we obtain 
  \begin{align}
  &\esssup_{t\in (0,T)} \intl_{B_{\lambda ^{ \frac{1}{4} k} }   }  |u_{k}(t)|^2  dx+  
  \intl_{0}^{T} \intl_{B_{ \lambda ^{ \frac{3}{5}k} }}  | \nabla u_{k}|^2 dxdt
  \cr
 &\qquad \qquad \qquad \le c K_0^2\Big(1+|\!|\!|b|\!|\!|^6 \max\{ T^{ \frac{13}{9}}, T\}\Big)\lambda ^{ - \frac{1}{20} k}.
\label{4.12}
\end{align}

\hspace{0.5cm}
Next, let $ l\in \N$ be fixed. Then \eqref{4.12} implies for all $ k \ge l$
\begin{align}
&\| u_k\|_{ L^\infty(0,T; L^2(B_{ \lambda^{ \frac{1}{4}l} }))} + \| \nabla u_k\|_{L^2(B_{ \lambda ^{ \frac{3}{5}l}}\times (0,T))}
\cr
&\qquad \qquad \le \| u_k\|_{ L^\infty(0,T; L^2(B_{ \lambda^{ \frac{1}{4}k} }))} + \| \nabla u_k\|_{L^2(B_{ \lambda ^{ \frac{3}{5}k}}\times (0,T))}
\cr
&\qquad \qquad \le c K_0\Big(1+|\!|\!|b|\!|\!|^3 \max\{ T^{ \frac{13}{18}}, T^{  \frac{1}{2} }\}\Big) \lambda ^{ - \frac{1}{40} k}.
\label{4.13}
\end{align}
Thus, by means of triangular inequality we find for each $ N\in \N, N >l$
\begin{align*}
&\Big\| \sum_{k=0}^{N} u_k\Big\|_{ L^\infty(0,T; L^2(B_{ \lambda^{ \frac{1}{4}l} }))}
+ \Big\|  \sum_{k=0}^{N} \nabla u_k \Big\|_{L^2(B_{ \lambda ^{ \frac{3}{5}l}}\times (0,T))}
\\
& \qquad \le \sum_{k=0}^{l-1}   \| u_k\|_{ L^\infty(0,T; L^2(\R^{3}))} + 
 \sum_{k=0}^{l-1}\| \nabla u_k\|_{L^2(\R^{3}\times (0,T))}
 \\
&\qquad \qquad + \sum_{k=l}^{N}   \| u_k\|_{ L^\infty(0,T; L^2(B_{ \lambda^{ \frac{1}{4}l} }))} + 
 \sum_{k=0}^{N}\| \nabla u_k\|_{L^2(B_{ \lambda ^{ \frac{3}{5}l}}\times (0,T))}
\\
 &\qquad \le c K_0 \lambda ^{ \frac{1}{2} l} +  c K_0\Big(1+ |\!|\!|b|\!|\!|^3 \max\{ T^{ \frac{13}{18}}, T^{  \frac{1}{2} }\}\Big)
 \\
 &\qquad \le c K_0\Big(\lambda ^{ \frac{1}{2} l}+ |\!|\!|b|\!|\!|^3 \max\{ T^{ \frac{13}{18}}, T^{  \frac{1}{2} }\}\Big). 
\end{align*}
Therefore,  $u^N= \sum_{k=0}^{N} u_k \rightarrow  u$ in $ V^2_{ \rm loc}(\R^{3}\times [0,T])$ as $ N \rightarrow\infty$. 
It is readily seen that $ u$ is a weak solution to \eqref{1.1}--\eqref{1.3},  and by virtue of the above  estimate we see that   for every $ 1 \le  \rho < \infty$
\begin{equation}
\| u\|_{ L^\infty(0, T; L^2(B_{\rho ^{  \frac{1}{4}}}))} + \| \nabla u\|_{ L^2(B_{ \rho^{\frac{3}{5} } }\times (0,T))}
\le c K_0\Big(\rho ^{ \frac{1}{2} }+ |\!|\!|b|\!|\!|^3 \max\{ T^{ \frac{13}{18}}, T^{  \frac{1}{2} }\}\Big).
\label{4.14}
\end{equation}
In particular, in \eqref{4.14} taking $ \rho =1$, and using  Sobolev's embedding theorem, we get   
\begin{equation}
\| u\|_{ L^4(0,T; L^3(B_1))} + \| u\|_{ V^2(B_1\times (0,T))} \le C_0 K_0\Big(1+ |\!|\!|b|\!|\!|^3 \max\{ T^{ \frac{13}{18}}, T^{  \frac{1}{2} }\}\Big)
\label{4.15}
\end{equation}
with a  constant $ C_0>0$ depending only on $ \lambda $.

\vspace{0.1cm}
\hspace{0.5cm}
It remains to show that $ u_\lambda = u$. Let $ N\in \N $, $ N \ge 4$. 
We set $ w^N = u^N - u^N_\lambda$. Recalling that $ b = b_\lambda$,  it follows that $ w^N$ solves  the system 
\begin{align}
\nabla \cdot  w^N &=0 \quad  \text{in}\quad  Q_{ \lambda ^{ -2} T},
\label{4.18}
\\
\partial _t w ^N + (b \cdot \nabla ) w^N -\Delta w^N  &= - \nabla \pi^N   \quad  \text{in}\quad  Q_{ \lambda ^{ -2} T},
\label{4.19}
\\
w^N &= w^N_0 \quad  \text{on}\quad  \R^{3}\times \{0\},
\label{4.20}
\end{align}

where
\begin{align*}
 w^N_0 &= \sum_{k=0}^{N} u_{0, k} -(u_{0, k })_\lambda 
 =\sum_{k=0}^{N} \PP( u_{0} \psi _k) -(\PP( u_{0} \psi _k))_\lambda  
\\
& = 
 u_{0} \sum_{k=0}^{N}\psi _k  -  \Big(u_{0} \sum_{k=0}^{N}\psi _k\Big)_{ \lambda }
 + \nabla \mathcal{N} \ast (u_0 \cdot \nabla \sum_{k=0}^{N}\psi _k) - 
\Big( \nabla \mathcal{N} \ast (u_0 \cdot \nabla \sum_{k=0}^{N}\psi _k) \Big)_\lambda
\\
& = u_{0} \Big(\sum_{k=0}^{N}\psi _k - \Big(\sum_{k=0}^{N}\psi _k\Big)(\lambda \cdot )\Big)   
 + \nabla \mathcal{N} \ast (u_0 \cdot \nabla \sum_{k=0}^{N}\psi _k) - 
\Big( \nabla \mathcal{N} \ast (u_0 \cdot \nabla \sum_{k=0}^{N}\psi _k) \Big)_\lambda,
\end{align*}
where $ \mathcal{N}= \frac{1}{4\pi | x|}$ stands for the Newton potential. For obtaining the third line in the above equalities we used the fact that $ (u_{ 0})_\lambda = u_0$. 
Owing to  $ \sum_{k=0}^{N}\psi _k= 1$ in $ B_{ \lambda ^{ N-3}}$ we have  
\begin{equation}
\Big(\sum_{k=0}^{N}\psi _k - \Big(\sum_{k=0}^{N}\psi _k\Big)(\lambda \cdot )\Big)=0\quad  \text{ in}\quad B_{ \lambda ^{ N-4}}.
\label{4.20a}
\end{equation}

\hspace{0.5cm}
Let $ \lambda ^{ \frac{3}{5}N } \le r < \rho \le \lambda ^{ \frac{3}{5}(N +1)} $ be arbitrarily chosen, but fixed.  Let $ \phi\in  C^{\infty}_{\rm c}(\R^{3})$ denote a cut off function such that $ 0 \le \phi \le 1$ in $ \R^{3}$, 
 $ \phi \equiv 1$ on $ B_{r}$, $ \phi \equiv 0$ in $ \R^{3}  \setminus B_{\rho }$, and $ | \nabla ^2 \phi |+ | \nabla \phi |^2  \le c (\rho -r)^{ -  2}$ in $ \R^{3}$. 
  Without loss of generality we may assume that $ \lambda ^{ \frac{3}{5}(N +1)} \le \lambda ^{ N-4}$. 
 Thus,  in view of \eqref{4.20a} we  infer that
   $ w_0^N$ is a gradient field in $ B_\rho $, and therefore 
  \begin{equation}
w_0^N- E_{ B_\rho }^{ \ast} (w_0^N)=0\quad  \text{ a.\,e. in $ B_\rho $}. 
 \label{4.20b}
 \end{equation}
 By a similar reasoning we have used to prove \eqref{4.10} we get the estimate  
\begin{align}
&  \| w^N\|^2_{ L^2(0,\lambda ^{ -2}T; L^6(B_r))}  +  \intl_{0}^{\lambda ^{ -2}T} \intl_{B_r}  | \nabla w^N|^2 dxdt
\cr
&\qquad \le cK_0^2 (1+|\!|\!|b|\!|\!|^6)\max\{ T^{ \frac{13}{9}}, T\}(\rho -r) ^{ - 6} \lambda ^{ \frac{17}{5}N}  
+ \frac{1}{2}\intl_{0}^{\lambda ^{ -2} T} \intl_{B_\rho }  | \nabla w^N|^2dxdt . 
\label{4.29}
\end{align}
Once more applying an iteration argument, together with the latter estimate, we deduce from \eqref{4.29}
\begin{equation}
  \| w^N\|^2_{ L^2(0,\lambda ^{ -2}T; L^6(B_{ \lambda ^{ \frac{3}{5} N}}))}  \le 
  c K_0^2   
(1+ |\!|\!|b|\!|\!|^6) \max\{ T^{ \frac{13}{9}}, T\}\lambda ^{ -\frac{1}{5} N}. 
\label{4.30a}
\end{equation}   
Accordingly, for all $ 0<\rho < \infty$, 
\[
 w^N \rightarrow 0  \quad  \text{{\it in}}\quad  L^2(0, \lambda ^{ -2}T; L^6(B_\rho ))\quad  \text{{\it as}}\quad  N \rightarrow +\infty. 
\]
On the other hand, observing that $ w^N = u^N - (u^N)_\lambda \rightarrow  u- u_\lambda  $ 
in $ L^2(0,\lambda ^{ -2}T; L^6(B_\rho ))$ as $ N \rightarrow \infty$, we conlude that $ u= u_\lambda $. 
This completes the proof of the  theorem.  \hfill \Beweisende 

\section{Proof of Theorem\,\ref{thm1.4}}
\label{sec:-6}
\setcounter{secnum}{\value{section} \setcounter{equation}{0}
\renewcommand{\theequation}{\mbox{\arabic{secnum}.\arabic{equation}}}}

We divide the proof in three steps. Firstly, given a $ \lambda $-DSS function  
$ b\in L^{ \frac{18}{5}}_{ loc}([0,\infty); L^3_{ \rm loc}(\R^{3}) )$ 
we get the existence of a unique $ \lambda $-DSS local  solution  with projected pressure $ u$ to the linearized system \eqref{4.1}--\eqref{4.3}, 
replacing $b $ by $ R_\var b$ therein (cf. appendix for the notion of the mollification $ R_\var $). 
Secondly,  based on the first step we may construct  a mapping $ \mathcal{T}: M \rightarrow  M$, which is continuous and compact. Application of  Schauder's fixed point theorem gives a local suitable solution with projected pressure to the approximated Navier-Stokes equation. Thirdly, 
letting $ \var \rightarrow 0^+$ in the weak formulation and in the local energy inequality \eqref{4.0b}, we obtain the existence of the desired local Leray solution with projected pressure to \eqref{1.1}--\eqref{1.3}.

\hspace{0.5cm}
We set 
\begin{equation}
T:=\min \Big\{\frac{1}{64C_0^{ 6}  K_0^6}, 
\Big(\frac{1}{64C_0^{ 6}  K_0^6}\Big)^{ \frac{9}{13}}\Big\}.
\label{6.0}
\end{equation}
Furthermore, set  $ X= L^3_{ \lambda -DSS}(Q)\cap  L^{ \frac{18}{5}}(0,T; L^3_{loc, \sigma}(\R^{3}))$  
equipped with the norm 
\[
|\!|\!| v|\!|\!| := \| v\|_{ L^{ \frac{18}{5}}(0,T; L^3(B_1)) },\quad  v\in X. 
\]
\vspace{0.5cm}  
Then we define,
\[
M= \Big\{ b\in X\,\Big|\,  |\!|\!| b|\!|\!| \le 2C_0 K_0\Big\}. 
\]
We now fix $ 0< \var < \lambda -1$.  For $ b\in M$  we set
\[
b_\var := R_{ \var }b, 
\]
where $ R_{ \var }$  stands for the mollification  operator defined in the appendix below. 
According to Theorem\,\ref{thm4.1} there exists a unique $ \lambda $-DSS  solution   $ u \in X$ to \eqref{4.1}--\eqref{4.3} with 
$ b _\var $ in place of $ b$.   Observing \eqref{4.15}, it follows that  
\begin{equation}
\| u \|_{ L^4(0,T; L^3(B_1))} + \| u \|_{ V^2(B_1\times (0,T))} \le C_0 K_0\Big(1+ |\!|\!|b|\!|\!|^3 \max\{ T^{ \frac{13}{18}}, T^{  \frac{1}{2} }\}\Big). 
\label{6.2}
\end{equation}
In view of \eqref{A.3a} having $ |\!|\!|b _\var |\!|\!|^3 \le \lambda ^{ \frac{5}{3}} |\!|\!|b |\!|\!|^3$, 
\eqref{6.2} together with \eqref{6.0} implies that 
\[
|\!|\!| u |\!|\!| \le 2 C_0 K_0,
\]  
and thus $ u \in M$. By setting $ \mathcal{T}_\var  (b):= u$  defines  a mapping   $ \mathcal{T}_\var : M \rightarrow  M$. 

\hspace{0.5cm}
{\it  $ \mathcal{T}_\var $ is closed}. In fact, let $\{ b_k\}$ be a sequence  in $ M$ such that $ b _k \rightarrow  b$ 
in $ X$ as $  k \rightarrow \infty$, and let $ u_k :=\mathcal{T}_\var (b_k)$, $ k\in \N$, such that $ u_k \rightarrow u$ in 
$ X$ as $ k \rightarrow \infty$.  From \eqref{6.2} it follows that $ \{ u_k\}$ is bounded in $ V^2_\sigma (B_1\times (0,T))$, 
and thus, eventually passing to a subsequence, we find that  $ u_k \rightarrow  u$ weakly in   $ V^2_\sigma (B_1\times (0,T))$ as $ k \rightarrow \infty$. Since $ u_k$ solves 
\eqref{4.1}--\eqref{4.3}  with $ b_{ k,\var }= R_{ \var } b_k$  in place of $ b$, from the above convergence properties we deduce that $ u\in M\cap V^2_\sigma (B_1\times (0,T))$ 
solves \eqref{4.1}--\eqref{4.3}. Accordingly, $ u= \mathcal{T}_\var (b)$. 

\hspace{0.5cm}
{\it  $ \mathcal{T}_\var (M) $ is relative compact in} $ X$.  To see this, let $ \{ u_k= \mathcal{T}_\var (b_k)\} \subset 
\mathcal{T}_\var (M)$ be any sequence.  Then $ u_k \in L^2_{ loc, \sigma}(\R^{3}\times [0,\infty))$ is a 
$ \lambda $-DSS local suitable weak    solution with projected pressure to 
\begin{align}
\nabla \cdot  u_k &=0 \quad  \text{in}\quad  Q,
\label{6.4}
\\
\partial _t u_k  + (b_{ k,\var } \cdot \nabla ) u_k -\Delta u_k  &= - \nabla \pi_k \quad  \text{in}\quad  Q,
\label{6.5}
\\
u_k  &= u_0 \quad  \text{on}\quad  \R^{3}\times \{0\}.
\label{6.6}
\end{align}
Introducing the local pressure, we have  
\begin{equation}
\partial _t v_k  + (b_{ k,\var } \cdot \nabla ) u_k -\Delta u_k  = - \nabla \pi_{ 1, k} - \nabla \pi_{2, k}  \quad  
\text{in}\quad  B_2\times (0,T),
\label{6.7}
\end{equation}
where $ v_k = u_k + \nabla p_{ h, k}$, and 
\begin{align*}
&\qquad \qquad \qquad \qquad  \nabla p_{ h,k} = - E^{ \ast}_{ B_2} (u_k), 
\\
&\nabla p_{ 1,k} = - E^{ \ast}_{ B_2} ( (b_{ k,\var } \cdot \nabla ) u_k),\quad  
\nabla p_{ 2,k} = E^{ \ast}_{ B_2} ( \Delta u_k). 
\end{align*}
Thus, \eqref{6.5} implies that $ v_k' = \nabla \cdot (- b_{ k,\var } \otimes u_k + \nabla u_k - p_{ 1,k}I- p_{ 2,k}I)$ in 
$ B_2\times (0,T)$.  Since $ b_k, u_k\in M$ we get the estimate 
\[
\| - b_{ k, \var } \otimes u_k + \nabla u_k - p_{ 1,k}I- p_{ 2,k}I\|_{ L^{ \frac{9}{5}}(0,T; L^{ \frac{3}{2}}(B_2))} 
\le c (1+ C_0^2 K_0^2).   
\]
Furthermore, by means of the reflexivity of $ L^2(0,T; W^{1,\, 2}(B_2))$, and using Banach-Alaoglu's theorem we get 
a subsequence $ \{ u_{ k_j}\}$ and a function $ u\in M \cap V^2_{loc, \sigma}(\R^{3}\times [0,T])$  such that 
\begin{align*}
 u_{ k_j} &\rightarrow u  \quad  \text{{\it weakly in}}\quad L^2(0,T; W^{1,\, 2}(B_2)),
\\
 u_{ k_j} &\rightarrow u  \quad  \text{{\it weakly$ ^{ \ast}$ in}}\quad  L^\infty(0,T; L^2(B_2))\quad  \text{{\it as}}\quad  j \rightarrow \infty.
\end{align*} 
In particular, we have for almost every  $ t\in (0,T)$
\begin{equation}
 u_{ k_j}(t) \rightarrow u(t)  \quad  \text{{\it weakly in}}\quad  L^2(B_2)\quad  \text{{\it as}}\quad  j \rightarrow \infty. 
\label{6.9}
\end{equation}

In addition, verifying that $ \{ v_{ k_j}\}$ is bounded in $ V^2(B_2\times (0,T))$, by Lions-Aubin's compactness lemma 
we see that 
\begin{equation}
 v_{ k_j} \rightarrow v  \quad  \text{{\it in}}\quad  L^2(B_2\times (0,T))\quad  \text{{\it as}}\quad  j \rightarrow +\infty, 
\label{6.10}
\end{equation}
where $ v= u+ \nabla p_{ h}$, and $ \nabla p_h = - E^{ \ast}(u)$.  Now, let $ t\in (0,T)$  be fixed such that \eqref{6.9} is satisfied. 
Then 
\begin{equation}
 \nabla p_{h,  k_j}(t) \rightarrow \nabla p_h (t)  \quad  \text{{\it weakly in}}\quad  L^2(B_2)\quad  \text{{\it as}}\quad  j \rightarrow \infty.  
\label{6.11}
\end{equation}
Since $ p_{ h,k}$ is harmonic in $ B_2$,  from \eqref{6.11} we deduce that 
\begin{equation}
 \nabla p_{h,  k_j}(t) \rightarrow \nabla p_h (t)  \quad  \text{{\it a.\,e. in}}\quad  B_2\quad  \text{{\it as}}\quad  j \rightarrow \infty.  
\label{6.12}
\end{equation} 
On the other hand, using the mean value property of harmonic functions, we see that $ \{ \nabla p_{h, k}\}$ 
is bounded in $ L^\infty(B_1\times (0,T))$. Appealing to Lebesgue's theorem of dominated convergence,  we infer from \eqref{6.12} 
that 
\begin{equation}
 \nabla p_{ h, k_j} \rightarrow \nabla p_h  \quad  \text{{\it in}}\quad  L^2(B_1\times (0,T))\quad  \text{{\it as}}\quad  j \rightarrow \infty.
\label{6.13}
\end{equation}   
Now combining \eqref{6.10} and \eqref{6.13}, we obtain $ u_{ k_j} \rightarrow  u$ in $ L^2(B_1\times (0,T))$. 
Recalling that $ \{ u_{ k_j}\}$ is bounded in $ V^2(B_1\times (0,T))$ , we get the desired convergence property 
$ u_{ k_j} \rightarrow u$ in $ X$ as $ j \rightarrow \infty$.  To see this we argue as follows. Eventually 
passing to a subsequence,  we may assume that $ u_{ k_j} \rightarrow  u$ almost everywhere in $ B_1\times (0,T)$. 
Let $ \var >0$ be arbitrarily chosen. We denote  $ A_m= \{(x,t)\in B_1\times (0,T) \,|\, \exists\,j \ge m: | u_{ k_j}(x,t) - u(x,t)|>\var \}$. Clearly, $ \cap_{ m=1}^\infty A_m $ is a set of Lebesgue measure zero. Thus $ \mes A_m \rightarrow 0$ as $ m
\rightarrow  \infty$. We now get the following estimate  
\begin{align*}
&\| u_{ k_j} - u\|_{ L^{ \frac{18}{5}}(0,T; L^3(B_1))} = 
\\
& \le  \| (u_{ k_j} - u)\chi _{ A_m}\|_{ L^{ \frac{18}{5}}(0,T; L^3(B_1))}+ \| (u_{ k_j} - u)\chi _{ A_m^c}\|_{ L^{ \frac{18}{5}}(0,T; L^3(B_1))}
\\
& \le \|u_{ k_j} - u \|_{ L^{ \frac{168}{45}}(0,T; L^{ \frac{28}{9}}(B_1))} \| \chi_{ A_m} \|_{ L^{ \frac{504}{5}}(0,T; L^{ 84}(B_1))}+  \| (u_{ k_j} - u)\chi _{ A_m^c}\|_{ L^{ \frac{18}{5}}(0,T; L^3(B_1))}. 
\\
& \le c (\mes A_m)^{ \frac{5}{504}} + c \var. 
\end{align*}
This shows that $ |\!|\!| u_{ k_j} - u|\!|\!| \rightarrow 0$ as $ j \rightarrow \infty$.  Applying Schauder's fixed point theorem, we get a function $ u_\var \in M$ such that $ u_\var  =\mathcal{T}_\var (u_\var )$. Thus, $ u_\var $ is a local suitable weak  solution with projected pressure to 
\begin{align}
\nabla \cdot  u_\var  &=0 \quad  \text{in}\quad  Q,
\label{6.15}
\\
\partial _t u_\var   + (R_\var u_\var  \cdot \nabla ) u_\var  -\Delta u_\var   &= - \nabla \pi_\var  \quad  \text{in}\quad  Q,
\label{6.16}
\\
u_\var   &= u_0 \quad  \text{on}\quad  \R^{3}\times \{0\}.
\label{6.17}
\end{align}

In particular, we have the a-priori estimate 
\begin{equation}
\| u_\var \|_{ L^4(0,T; L^3(B_1))} + \| u_\var \|_{ V^2(B_1\times (0,T))} \le 2C_0K_0. 
\label{6.21}
\end{equation}

Let $ \{ \var _j\}$ be a sequence of positive numbers in $ (0, \lambda -1)$. 
Since $ u_\var $ is $ \lambda $-DSS we have $ u_{ \var }(x,t)= u_{ \var , \lambda }(x,t)$ for almost every $ (x,t)\in Q$. 
Thus, there exists a set of measure zero $ S \subset (0, +\infty)$  such that for all $ t\in [0,+\infty)  \setminus S$  
\[
u_{ \var_j }(x,t)= u_{ \var_j , \lambda }(x,t) = \lambda^{ k} u_{ \var _j}(\lambda^{ k} x, \lambda ^{ 2k}t)\quad \text{ for a.\,e. $ x\in \R^{3}$},\quad  \forall\,k\in \Z, \forall\,j\in \N. 
\] 

\hspace{0.5cm}
Clearly, $ t\in (0,+\infty)  \setminus S$ iff $ \lambda ^2 t\in (0,+\infty)  \setminus S$.  Indeed, let $ t\in N^c $. Then 
$\lambda  u_{ \var _j} (\lambda x, t) = u_{ \var _j}( \lambda  ^2x, \lambda ^2 t) $ for almost every $ x\in \R^{3}$. By means of the reflexivity we get a sequence $ \var _j \rightarrow 0^+ $ as $ j \rightarrow \infty$ and $ u\in V^2_{ loc, \sigma}(\R^{3}\times [0,T])$ such that 
\begin{align*}
 u_{ \var _j } &\rightarrow u \quad  \text{{\it weakly in}}\quad  L^2(0,T; W^{1,\, 2}(B_1))\quad  
 \text{{\it as}}\quad  j \rightarrow +\infty,
\\
 u_{ \var _j} &\rightarrow u  \quad  \text{{\it weakly$^{ \ast} $ in}}\quad  L^\infty(0,T; L^2(B_1))\quad  \text{{\it as}}\quad  j \rightarrow +\infty.
\end{align*}
Arguing as in the proof the compactness of $ \mathcal{T}_\var $, we infer 
\[
u_{ \var _j}  \rightarrow  u\quad  \text{{\it  in}}\quad L^{ \frac{18}{5}}(0, T; L^3(B_1))\quad  
\text{{\it as}}\quad  j \rightarrow 0^+.   
\]
Note that $ u$ is DSS, since $ u$ is obtained as a limit of sequence DSS functions. 

\hspace{0.5cm}
Together with  Lemma\,\ref{lemA.3} we see that   
\begin{equation}
R_{ \var _j} u_{ \var _j}  \rightarrow  u\quad  \text{{\it  in}}\quad L^{ \frac{18}{5}}(0, T; L^3(B_1))\quad  
\text{{\it as}}\quad  j \rightarrow 0^+.   
\label{6.20}
\end{equation}
This shows that  $ u\in L^2_{ loc,  \sigma}(\R^{3}\times [0,+\infty))$ is a local Leray solution with projected pressure to \eqref{1.1}--\eqref{1.3}.  \hfill \Beweisende


\hspace{0.5cm}
$$\mbox{\bf Acknowledgements}$$
Chae was partially supported by NRF grants 2016R1A2B3011647, while Wolf has been supported 
supported by the German Research Foundation (DFG) through the project WO1988/1-1; 612414.

\appendix
\section{Mollification for DSS functions}
\label{sec:-A}
\setcounter{secnum}{\value{section} \setcounter{equation}{0}
\renewcommand{\theequation}{\mbox{A.\arabic{equation}}}}

Let $1< \lambda < +\infty$. Let $ u \in L^s_{ \lambda - DSS}(\R^{3})$. Let $ \rho \in C^{\infty}_{\rm c}(B_1) $ denote the standard mollifying kernel such that $  \intl_{\R^{3}}\rho dx =1$. 
For $ 0<\var <\lambda -1$ we define 
\[
(R_\var u)(x,t) = \frac{1}{(\,\sqrt{t} \var)^3}
\intl_{B_{ \,\sqrt{t}\var }} u(x-y,t) \rho \Big( \frac{y}{ \,\sqrt{t} \var }\Big) dy,\quad  (x,t)\in Q.  
\]

We have the following 

\begin{lem}
\label{lemA.1}
$ R_\var $ defines a bounded operator from $ L^s_{ \lambda - DSS}(Q)$ into itself. Furthermore, 
 for all $ u\in L^s_{ \lambda - DSS}(Q)$  it holds for all $ (x.t)\in Q$
 \begin{equation}
 | (R_\var u)(x,t)| \le  c\{ \,\sqrt{t} \var\} ^{ - \frac{3}{s}}  \|u(\cdot , t)  \|_{ L^s(B_{ \,\sqrt{t}\var }(x))}
 \label{A.1}
 \end{equation}
 with an constant $ c>0$ depending on $ s$ only. 
\end{lem}

{\bf Proof}: Let $ u\in L^s_{ \lambda - DSS}(Q)$. First we will verify that $ R_\var u$ is $ \lambda $-DSS. 
Indeed, using the transformation formula of the Lebesgue integral, we calculate  for any $ (x,t)\in Q$, 
\begin{align*}
\lambda (R_\var  u)(\lambda x, \lambda ^2t) &=  \frac{1}{\lambda ^2(\,\sqrt{t} \var)^3}
\intl_{B_{\lambda  \,\sqrt{t}\var }} u(\lambda x-y, \lambda ^2t) \rho \Big( \frac{y}{ \lambda \,\sqrt{t} \var }\Big) dy,
\\
&= \frac{1}{(\,\sqrt{t} \var)^3}\intl_{\R^{3}} \lambda u(\lambda (x-y), \lambda ^2 t) \rho \Big( \frac{y}{  \,\sqrt{t} \var }\Big) dy  
\\
\\
&= \frac{1}{(\,\sqrt{t} \var)^3}\intl_{\R^{3}} u(x-y,  t) \rho \Big( \frac{y}{  \,\sqrt{t} \var }\Big) dy =
(R_\var  u) (x,t). 
\end{align*}

\hspace{0.5cm}
Firstly, let  $ \lambda ^{ -2}< t \le  1$.  
Noting that $ (R_\var  u)(\cdot , t) = u(\cdot ,t) \ast \rho _{ \,\sqrt{t} \var }$, 
where $ \rho _{ \,\sqrt{t} \var }(y) = \frac{1}{ (\,\sqrt{t} \var )^3 } \rho \Big( \frac{y}{ \,\sqrt{t} \var }\Big)$, 
recalling that $ \var < \lambda -1$,  by means of Young's inequality we find 
\[
\| (R_\var  u)(\cdot , t)\|^s_{ L^s(B_1)} \le  \| u(\cdot , t)\|^s_{ L^s(B_{ 1+\var })} \| \rho _{ \,\sqrt{t} \var } \|^s_{ L^1} 
= \| u(\cdot , t)\|^s_{ L^s(B_{\lambda })}. 
\]
 Integrating the above inequality over $ (\lambda ^{ -2},1)$, and using a suitable change of coordinates,  we obtain  
\begin{align*}
\| R_\var  u\|_{ L^s(B_1 \times (\lambda ^{ -2}, 1))} &\le \| u\|_{ L^s(B_\lambda  \times (\lambda ^{ -2}, 1))} 
\\
&= 
\| u\|_{ L^s(B_1 \times (\lambda ^{ -2}, 1) )} + \| u\|_{ L^s(B_\lambda   \setminus B_{1 }\times (\lambda ^{ -2}, 1) )} 
\\
&=
\| u\|_{ L^s(B_1 \times (\lambda ^{ -2}, 1) )} +\lambda ^{ \frac{5-s}{s}}\| u\|_{ L^s(B_1   \setminus B_{\lambda ^{ -1} }\times (\lambda ^{ -4}, \lambda ^{ -2}) )}. 
\end{align*}

\hspace{0.5cm}
Secondly, for $ 0< t < \lambda ^{ -2}$ we estimate 
\[
\| (R_\var  u)(\cdot , t)\|^s_{ L^s(B_1  \setminus B_{ \lambda ^{ -1}})} \le  \| u(\cdot , t)\|^s_{ L^s(B_{ \lambda }  \setminus B_{ \lambda ^{ -1}})} \| \rho _{ \,\sqrt{t} \var } \|^s_{ L^1} 
= \| u(\cdot , t)\|^s_{ L^s(B_{ \lambda }  \setminus B_{ \lambda ^{ -1}})}. 
\]
Integration over $ (0, \lambda ^{ -2})$ in time yields 
\begin{align*}
\| R_\var  u\|_{ L^s(B_1  \setminus B_{ \lambda ^{ -1}}\times (0,\lambda ^{ -2}))} &\le \| u\|_{ L^s(B_\lambda  \setminus B_{ \lambda ^{ -1}}  \times (0, \lambda ^{ -2}))} 
\\
&= 
\| u\|_{ L^s(B_1  \setminus B_{ \lambda ^{ -1}} \times (0, \lambda ^{ -2}) )} + 
\| u\|_{ L^s(B_\lambda   \setminus B_{1 }\times (0, \lambda ^{ -2}) )} 
\\
&=
\| u\|_{ L^s(B_1  \setminus B_{ \lambda ^{ -1}} \times (0, \lambda ^{ -2}) )} +\lambda ^{ \frac{5-s}{s}}\| u\|_{ L^s(B_1   \setminus B_{\lambda ^{ -1} }\times (0, \lambda ^{ -4})}. 
\end{align*}
Combining the last two estimates, we get
\[
\|R_\var  u \|_{ L^s(Q_1  \setminus Q_{ \lambda ^{ -1}})} \le  (1+ \lambda ^{ \frac{5-s}{s}}) \| u\|_{ L^s(Q_1  \setminus Q_{ \lambda ^{ -1}})}.  
\] 
This shows that $ R_\var : L^s_{ \lambda -DSS}(Q) \rightarrow  L^s_{ \lambda -DSS}(Q)$ is bounded. 

\hspace{0.5cm}
The inequality \eqref{A.1} follows immediately from the definition of $ R_\var u$ with the help of  H\"older's inequality.   
 \hfill \Beweisende 

\begin{rem}
\label{remA.2} 
Arguing as in the proof of Lemma\,\ref{lemA.1},  we get for any 
$ u\in L^3_{ \lambda -DSS}(Q)\cap L^{ \frac{18}{5}}(0,T; L^3(B_1))$, $ 0<T<1$
\begin{equation}
\| R_\var  u\|_{ L^{ \frac{18}{5}}(0,T; L^3(B_1))} \le  \lambda^{ \frac{5}{9}} 
\| u\|_{ L^{ \frac{18}{5}}(0,T; L^3(B_1))}. 
\label{A.3a}
\end{equation}
\end{rem}

\begin{lem}
\label{lemA.3}
Let $ u\in L^3_{ \lambda -DSS}(Q)\cap L^{ \frac{18}{5}}(0,T; L^3(B_1))$, $ 0<T \le 1$. Then 
\begin{equation}
 R_\var u  \rightarrow u  \quad  \text{{\it in}}\quad  L^{ \frac{18}{5}}(0,T; L^3(B_1))\quad  \text{{\it as}}\quad  \var  \rightarrow 0^+. 
\label{A.2}
\end{equation} 

\end{lem}

{\bf Proof}:  First by the absolutely continuity of the Lebesgue integral we see that for almost all $ t\in (0,T)$
\[
 (R_\var u)(\cdot , t) \rightarrow u(\cdot , t) \quad  \text{{\it in}}\quad  L^3(B_1)\quad  \text{{\it as}}\quad  \var  \rightarrow 0^+.  
\]

Let $ A \subset (0,T)$ be any Lebesgue measurable set. By Young's inequality of convolutions we get 
for almost all $ t\in (0,T)$
\[
 \intl_{A} \| (R_\var u)(\cdot, t)\|_{ L^3(B_1)}^{ \frac{18}{5}}    dt 
 \le \intl_{A} \| u(\cdot, t)\|_{ L^3(B_\lambda )}^{ \frac{18}{5}}    dt
\]
Since $ u \in L^{ \frac{18}{5}}(0,T; L ^3(B_\lambda ))$,  the assertion \eqref{A.2}  follows by the aid of Vitali's convergence 
lemma.  \hfill \Beweisende 

\section{Weak trace for time dependent  $ \lambda $-DSS functions}
\label{sec:-B}
\setcounter{secnum}{\value{section} \setcounter{equation}{0}
\renewcommand{\theequation}{\mbox{B.\arabic{equation}}}}

Let $ 1<\lambda < +\infty$. A measurable function   $ u: Q \rightarrow \R^3$ is said to be $ \lambda $-DSS,  if 
for almost every $ (x,t)\in Q$ 
\begin{equation}
u(x,t) = \lambda u(\lambda x,\lambda ^2t).  
\label{B.1}
\end{equation}.  
We denote by  $ M(u)$  the set of all $ t\in [0, +\infty)$ such that for all $ k\in \Z$
\begin{equation}
u(x,t) = \lambda^{ k} u(\lambda^{ k} x,\lambda^{ 2k}t)\quad  \text{ for a.\,e. \, $ x\in \R^{3}$}.  
\label{B.2}
\end{equation}.  

\begin{lem}
\label{lemB.1}
The set $ [0,+\infty)  \setminus M(u)$ is a set of Lebesgue measure zero.
\end{lem}

{\bf Proof}: For $ m\in N$ and $ k\in \N$ by $ A_{ m, k}$ we denote the set of all $ t\in [0, +\infty)$ such that 
\[
\mes \Big\{x\in \R^{3} \,\Big|\, u(x,t) = \lambda^{ k} u(\lambda^{ k} x,\lambda^{ 2k}t)\Big\} \ge \frac{1}{m}. 
\]
Since $ u$ is discretely self-similar,  we must have $ \mes (A_{ m,k})=0$. Since $ M(u)  \setminus [0, +\infty) = 
\cup_{ k\in \Z}\cup_{ m=1}^\infty A_{ m,k} $ the assertion follows. 
 \hfill \Beweisende

\begin{lem}
\label{lemB.2}
For every $ t\in [0,+\infty)$ it holds $ t\in M(u)$ iff $ \lambda ^{ 2} t \in M(u)$. 
\end{lem}

{\bf Proof}: Let $ t\in M(u)$. There exists a set $ P \subset \R^{3}$ with $ \mes (\R^{3}  \setminus P)=0$ such that 
\eqref{B.2} holds for all $ x\in P$. Define $ P_k = \{y=\lambda^{ k} x  \,|\, x \in P\}$, $ k\in \Z$. Cleary, 
$ \mes (\R^{3}  \setminus \cap_{ k\in \Z} P_k )=0$. Let  $ x\in \cap_{ k\in \Z} P_k$. Then $x,  \lambda ^{ -1}x\in P$, and therefore for all $ k\in \Z$ we get $u(\lambda ^{ -1} x, t) = \lambda u( x, \lambda ^2 t) = \lambda^{k+1} u(\lambda^k x, \lambda^{ 2+2k} t)$,  which is equivalent to 
\[
u(x, t) = \lambda ^{ k} u(\lambda ^kx, \lambda ^{ 2k} \lambda ^{ 2} t). 
\]
This shows that $ \lambda ^2t \in M(u)$. Similarly, we get the opposite direction.  \hfill \Beweisende     

\vspace{0.5cm}  
As an immediate consequence of Lemma\,\ref{lemB.1} we see that
\begin{equation}
t\in M(u) \quad  \Longleftrightarrow \quad  \quad \lambda ^{ 2k} t\in M(u)\quad  \forall\,k\in \Z.
\label{B.3}
\end{equation}

\hspace{0.5cm}
Let $ \{ v_j\}$ be a sequence 
in $ L^2_{ loc}(\R^{3})$. We say 
\[
 v_j \rightarrow v  \quad  \text{{\it weakly in}}\quad  L^2_{ loc}(\R^{3})\quad  \text{{\it as}}\quad  j \rightarrow +\infty 
\]
if for every $ 0<R<+\infty$
\[
 v_j \rightarrow v  \quad  \text{{\it weakly in}}\quad  L^2(B_R)\quad  \text{{\it as}}\quad  j \rightarrow +\infty. 
\]

\begin{lem}
\label{lemB.3}
Let $ \{ v_j\}$ be a sequence in $ L^2_{ loc}(\R^{3})$ such that for all $ 0<R<+\infty$
\begin{equation}
\sup_{ j\in \N} \| v_j\|_{ L^2(B_R)} <+\infty. 
\label{B.5}
\end{equation}
Then there exists a subsequence $ \{ v_{ j_m}\}$ and $ v\in L^2_{ loc}(\R^{3})$ such that 
\[
 v_{ j_m} \rightarrow v  \quad  \text{{\it weakly in}}\quad  L^2_{ loc}(\R^{3})\quad  \text{{\it as}}\quad  m \rightarrow +\infty. 
\]

\end{lem}

{\bf Proof}:  By induction and the reflexivity of $ L^2(B_m)$ we construct  a sequence of subsequences $ \{ v_{ j^{ (m)}_k}\} \subset \{ v_{ j^{ (m-1)}_k}\}$ and $ \{ v_{ j^{0}_k}\}=\{ v_j\}$ such that for some $ v_m \in L^2(B_k)$ it holds 
\[
 v_{ j^{ (m)}_k} \rightarrow v_m  \quad  \text{{\it in}} \quad  L^2(B_m)\quad  \text{{\it as}}\quad  k \rightarrow +\infty
\]
$ (m\in \N)$. Clearly, $ v_{ m}|_{ B_{ m-1}}= v_{ m-1} $.  This allows us to define $ v: \R^{3} \rightarrow \R$ be setting 
$ v=v_m$ on $ B_m$. Then by Cantor's diagonalization principle the subsequence   $ v_{ j_m}= v_{ j^{ (m)}_m}$ meets the requirements.   \hfill \Beweisende 

\vspace{0.2cm}
\hspace{0.5cm}
We denote $ \mathcal{V}=  L^\infty_{ loc}([0,+\infty); L_{ loc}^2(\R^{3}))$ the space of all measurable functions $ u: Q
\rightarrow \R$ such that $ u\in L^\infty(0, R^2, L^2(B_R))$ for all $0<R<+\infty$. By $ \mathcal{V}_{ \lambda -DSS}$
 we denote the space of all $ \lambda $-DSS functions $ u\in \mathcal{V}$.

\begin{lem}
\label{lemB.4}
Let $ u\in \mathcal{V}_{ \lambda -DSS}$. We assume that $ \| u(t) \|_{ L^2(B_R)} \le \| u\|_{ L^\infty(0,R ^2; L^2(B_R))}$
for all $ t\in (0,R^2)$, $ 0<R<+\infty$.  There exists a constant $ C>0$  such that  for every $ t\in  M$ 
\begin{equation}
\| u(t)\|^2_{ L^2(B_{R})} \le C\max
\Big\{R \| u\|_{ L ^\infty(0,1; L^2(B_1))}, \| u(t)\|_{  L^2(B_{ \,\sqrt{t}})}\Big\}.
\label{B.6}
\end{equation}

 \end{lem}

{\bf Proof}:  Let $ t\in M(u)$. Let $ k\in \Z$. Then by means of the transformation formula we get 
\begin{align*}
\intl_{A_k} | u(x,t)|^2 dx &= \lambda^{ 3k} \intl_{A_1} | u(\lambda^k x, t)|^2 dx 
= \lambda^{ k} \intl_{A_1} | \lambda ^{ k}u(\lambda^k x, \lambda^{ 2k} \lambda ^{ -2k}t)|^2 dx 
\\
&= \lambda^{ k} \intl_{A_1} | u( x,  \lambda ^{ -2k}t)|^2 dx 
\end{align*}
In case $ \lambda ^{ 2k} \ge t $ we get 
\[
\| u(t)\|^2_{ L^2(A_k)} \le \lambda ^k \| u\|_{ L ^\infty(0,1; L^2(B_1))}.  
\]
 On the contrary, if $  \lambda ^{ 2k}< t $ we find 
\[
\| u(t)\|^2_{ L^2(B_{ \lambda ^k})} \le \| u(t)\|_{  L^2(B_{ \,\sqrt{t}})}.  
\] 
Accordingly, 
 \[
\| u(t)\|^2_{ L^2(B_{ \lambda ^k})} \le c  \max
\Big\{\lambda ^k \| u\|_{ L ^\infty(0,1; L^2(B_1))},  \| u(t)\|_{  L^2(B_{ \,\sqrt{t}})}\Big\}.
\]
This yields \eqref{B.6}.  \hfill \Beweisende  

\begin{lem}
\label{lemB5}
Let $ u\in \mathcal{V}_{ \lambda -DSS}$. Furthermore, let $ F_{ ij}, g_i: Q \rightarrow \R$ such that 
$ F_{ ij}, g_{ i}\in L^1(Q_R)$ and for all $ 0<R<+\infty$, i, j=1,2,3. We suppose for all $ t\in [0,+\infty)$ the function  $u(\cdot , t)\in L^2_{ loc}(\R^{3}) $ 
with $ \nabla \cdot u(\cdot , t)=0$ in the sense of distributions, and 
that for all $ \varphi \in C^{\infty}_{\rm c}(Q) $  with $ \nabla \cdot \varphi =0$ the following identity holds true 
\begin{equation}
\intl_{Q}  u \frac{\partial \varphi }{\partial t} dxdt = \intl_{Q}   F: \nabla \varphi  + g\cdot \varphi dxdt. 
\label{B.7}
\end{equation}
Then, eventually redefining  $ u(t)$ for $ t$ in a set of measure zero, we have 
\begin{align}
& u \in C_w([0, +\infty); L^2(B_R))\quad  \forall\,0<R<+\infty,
\label{B.8}
\\
& M(u) = [0,+\infty). 
\label{B.9}
\end{align}

\end{lem}

{\bf Proof}: By $ L(u) \subset [0,+\infty)$ we denote the set of all Lebesgue points of $ u$, more precisely, we say 
$ t \in L(u)$,  if for every $ 0< R< +\infty$  
\[
\frac{1}{\var }  \intl_{t}^{t+\var }   u(\cdot , \tau ) d \tau \rightarrow u(\cdot ,t)\quad  \text{ {\it in}}\quad  L^2(B_R)\quad  \text{ {\it as}}\quad \var \rightarrow +\infty. 
\]
By  Lebesgue's differentiation theorem we have $ \mes ([0,+\infty)  \setminus L(u))=0$.  
Let $ t\in L(u)$. By a standard approximation argument 
we deduce from \eqref{B.7} that for every $ \varphi \in C^{\infty}_{\rm c}(Q)$ with $ \nabla \cdot u=0$ 
\begin{equation}
-\intl_{\R^{3}} u(t) \varphi (t)dx +    \intl_{0}^{t}  \intl_{\R^{3}} u \frac{\partial \varphi }{\partial t} dxds =   \intl_{0}^{t}  \intl_{\R^{3}}   F: \nabla \varphi  + g\cdot \varphi dxds. 
\label{B.10}
\end{equation}
Next, let $ \{ t_j\}$ be a sequence in $ M(u)\cap L(u)$ such that $ t_j \rightarrow t\in L(u)$ as $ j \in +\infty$.  Thank's to Lemma\,\ref{lemB.3} there exists a subsequence $ \{ t_{ j_m}\}$ and $ v\in L^2_{ loc}(\R^{3})$ such that 
\[
 u(t_{ j_m}) \rightarrow v  \quad  \text{{\it weakly in}}\quad  L^2_{ loc}(\R^{3})\quad  \text{{\it as}}\quad  m \rightarrow +\infty.  
\]
Thus, \eqref{B.10} implies for all $ \varphi \in C^{\infty}_{\rm c}(Q)$ with $ \nabla \cdot \varphi =0$
\begin{equation}
-\intl_{\R^{3}} v \varphi (t)dx +    \intl_{0}^{t}  \intl_{\R^{3}} u \frac{\partial \varphi }{\partial t} dxds =   \intl_{0}^{t}  \intl_{\R^{3}}   F: \nabla \varphi  + g\cdot \varphi dxds. 
\label{B.11}
\end{equation} 
On the other hand, recalling that $ t\in L(u)$,  we have the same identity as \eqref{B.11} replacing $ v$ by $ u(t)$ therein. This shows that 
for all $ \psi \in C^{\infty}_{\rm c, \sigma }(\R^{3})$
\[
\intl_{\R^{3}} (v- u(t)) \cdot \psi dx =0. 
\]
Consequently, $ v- u(t)$ is a harmonic function. On the other hand, by the lower semi continuity of the $ L^2$ norm we obtain from 
\eqref{B.6} that 
\begin{equation}
\| u(t)- v\|^2_{ L^2(B_{R})} \le C\max
\Big\{R \| u\|_{ L ^\infty(0,1; L^2(B_1))}, \| u(t)\|_{  L^2(B_{ \,\sqrt{t}})}\Big\}.
\label{B.12}
\end{equation} 
 Whence, $ v=u(t)$. In particular, $ u(s) \rightarrow u(t)$ weakly in $ L^2_{ loc}(\R^{3})$ as $ s\in M(u)\cap L(u) \rightarrow t$.  

\hspace{0.5cm}
Let $ t\in [0,+\infty)$. There exists a sequence $ \{ t_j\}$ in $ M(u)\cap L(u)$ such that $ t_j \rightarrow t$ as $ j \rightarrow +\infty$. 
 Thank's to Lemma\,\ref{lemB.3} there exists a subsequence $ \{ t_{ j_m}\}$ and $ v\in L^2_{ loc}(\R^{3})$ such that 
\[
 u(t_{ j_m}) \rightarrow v  \quad  \text{{\it weakly in}}\quad  L^2_{ loc}(\R^{3})\quad  \text{{\it as}}
 \quad  m \rightarrow +\infty.  
\]
Observing  \eqref{B.10} with $ t_{j_m}$  in place of $ t$ and letting $  m \rightarrow +\infty$, we obtain  for all $ \varphi \in C^{\infty}_{\rm c}(Q)$ with $ \nabla \cdot \varphi =0$
\begin{equation}
-\intl_{\R^{3}} v \varphi (t)dx +    \intl_{0}^{t}  \intl_{\R^{3}} u \frac{\partial \varphi }{\partial t} dxds =   \intl_{0}^{t}  \intl_{\R^{3}}   F: \nabla \varphi  + g\cdot \varphi dxds. 
\label{B.14}
\end{equation} 
On the other hand, by the lower semi continuity of the $ L^2$ norm  from 
\eqref{B.6} it follows that  
\begin{equation}
\| v\|^2_{ L^2(B_{R})} \le C\max
\Big\{R \| u\|_{ L ^\infty(0,1; L^2(B_1))}, \| u(t)\|_{  L^2(B_{ \,\sqrt{t}})}\Big\}.
\label{B.15}
\end{equation} 
For a second subsequence $ \{ t_{ j_m}'\}$ with limit $ w\in L^2_{ loc}(\R^{3})$ we derive the same property as $ v$ which leads to the fact that for all $ \psi \in C^{\infty}_{\rm c, \sigma }(\R^{3})$
\[
\intl_{\R^{3}} (v-w) \cdot \psi dx =0. 
\]
Consequently, $ v-w$ is a harmonic function. Now taking into account the estimate \eqref{B.12},  which is satisfied for $ w$ too, 
we infer $ v=w$.  Thus, the limit is uniquely determined. In case $ t\notin M(u)\cap L(u)$ we set $ u(t)=v$.   In particular, by the 
lower semi continuity of the norm we have for all $ t\in [0,+\infty)$ the estimate 
\begin{equation}
\| u(t)\|^2_{ L^2(B_{R})} \le C\max
\Big\{R \| u\|_{ L ^\infty(0,1; L^2(B_1))}, \| u(t)\|_{  L^2(B_{ \,\sqrt{t}})}\Big\}.
\label{B.16}
\end{equation} 
Furthermore, observing \eqref{B.11},  it follows that for all $ t\in [0,+\infty)$ and for all $ \varphi \in C^{\infty}_{\rm c}(Q)$ with $ \nabla \varphi =0$
  \begin{equation}
-\intl_{\R^{3}} u(t) \varphi (t)dx +    \intl_{0}^{t}  \intl_{\R^{3}} u \frac{\partial \varphi }{\partial t} dxds =   \intl_{0}^{t}  \intl_{\R^{3}}   F: \nabla \varphi  + g\cdot \varphi dxds. 
\label{B.17}
\end{equation} 

\hspace{0.5cm}
Next, let $ t\in [0,+\infty)$, and let $ \{ t_j\}$ be any sequence in $ [0,+\infty)$ with $ t_j \rightarrow t$ as $ j \rightarrow +\infty$. 
Once more applying Lemma\,\ref{lemB.3}, we get a subsequence $ \{ t_{ j_m}\}$ and $ w\in L^2_{ loc}(\R^{3})$ such that 
\[
 u(t_{ j_m}) \rightarrow w  \quad  \text{{\it weakly in}}\quad  L^2_{ loc}(\R^{3})\quad  \text{{\it as}}
 \quad  m \rightarrow +\infty.  
\]
Observing \eqref{B.17} with $ t_{ j_m}$ in place of $ t$ and letting $ m \rightarrow +\infty$,  it follows that 
  \begin{equation}
-\intl_{\R^{3}} w \varphi (t)dx +    \intl_{0}^{t}  \intl_{\R^{3}} u \frac{\partial \varphi }{\partial t} dxds =   \intl_{0}^{t}  \intl_{\R^{3}}   F: \nabla \varphi  + g\cdot \varphi dxds. 
\label{B.18}
\end{equation}   
Combining \eqref{B.18} and \eqref{B.17} and verifying \eqref{B.15} for $ w$ by a similar reasoning as above, we conclude $ w=u(t)$. 
This shows that $ u\in C_w([0,+\infty); L ^2_{ loc}(\R^{3}))$.    

\hspace{0.5cm}
It only remains to prove that $ M(u)= [0,+\infty)$. To see this let $ \{ t_j\}$ be a sequence in $ M(u)$ such that $ t_j \rightarrow t$. 
By using the transformation formula of the Lebesgue integral together with Lemma\,\ref{lemB.2} (cf. also \eqref{B.3} ),   we calculate for all $ \psi \in C^{\infty}_{\rm c}(\R^{3})$
\begin{align*}
\intl_{\R^{3}} u(x, t) \psi(x)  dx &= \lim_{ j\to \infty} \intl_{\R^{3}} u(x, t_j) \psi(x)  dx
\\
&= \lambda ^{ -3k}\lim_{ j\to \infty} \intl_{\R^{3}}  u(\lambda^{ -k} x, t_j) \psi(\lambda x)  dx
\\
&= \lambda ^{ -2k}\lim_{ j\to \infty} \intl_{\R^{3}}  u(x, \lambda ^{2k}t_j) \psi(\lambda x)  dx
\\
&=  \lambda ^{ -2k}\intl_{\R^{3}}  u(x, \lambda ^{2k}t) \psi(\lambda x)  dx
= \intl_{\R^{3}} \lambda ^{k} u(\lambda^k x, \lambda ^{ 2k}t) \psi( x)  dx.
\end{align*}
 This yields $u(x, t) =\lambda ^{k} u(\lambda^k x, \lambda ^{ 2k}t)  $ for almost every $ (x,t)\in Q$, and thus $ t\in M(u)$. 
  \hfill \Beweisende


\begin{thebibliography}{10}

\bibitem{BradTsai2015}
{\sc Z.~Bradshaw and T.-P. Tsai}, {\em Forward discretely self-similar
  solutions of the {N}avier-{S}tokes equations II}, arXiv:1510.07504v1 (2015), to appear in 
  Ann. I. H. Poincar\'e-AN.

\bibitem{CaKoNi1982}
{\sc L.~Caffarelli, R.~Kohn, and L.~Nirenberg}, {\em Partial regularity of
  suitable weak solutions of the {N}avier-{S}tokes equations}, Comm. Pure Appl.
  Math., 35 (1982), pp.~771--831.

\bibitem{GaSiSo1994}
{\sc G.~Galdi, C.~Simader, and H.~Sohr}, {\em On the stokes problem in
  lipschitz domains}, Annali di Mat. pura ed appl. (IV), 167 (1994),
  pp.~147--163.

\bibitem{JiaSve2014}
{\sc H.~Jia and V.~\v{S}ver\'{a}k}, {\em Local-in-space estimates near initial
  time for weak solutions of the {N}avier-{S}tokes equations and forward
  self-similar solutions}, Invent. Math., 196 (2014), pp.~233--265.

\bibitem{KiSe2007}
{\sc N.~Kikuchi and G.~A. Seregin}, {\em Weak solutions to the cauchy problem
  for the {N}avier-{S}tokes equations satisfying the local energy inequality},
  in Nonlinear equations and spectral theory, A.~M.~S. Transl., ed., no.~220 in
  Ser 2., Providence, RI, 2007, Amer. Math. Soc., pp.~141--164.

\bibitem{KoTa2001}
{\sc H.~Koch and D.~Tataru}, {\em Well-posedness for the {N}avier-{S}tokes
  equations}, Adv. Math., 157 (2001), pp.~22--35.

\bibitem{LemRie2002}
{\sc P.~G. Lemari{\`e}-Rieusset}, {\em Recent developments in the
  {N}avier-{S}tokes problem}, vol.~431, Chapman Hall/CRC, Boca Raton, FL, 2002.

\bibitem{Leray1934}
{\sc J.~Leray}, {\em Sur le mouvement d'un liquide visqueux emplissant
  l'espace}, Acta Math., 63 (1934), pp.~193--284.

\bibitem{Scheffer1976}
{\sc V.~Scheffer}, {\em Partial regularity of solutions to the
  {N}avier-{S}tokes equations}, Pacific J. Math., 66 (1976), pp.~535--552.

\bibitem{Tsai2014}
{\sc T.-P. Tsai}, {\em Forward discretely self-similar solutions of the
  {N}avier-{S}tokes equations}, Comm. Math. Phys., 328 (2014), pp.~29--44.


\bibitem{Wolf2015c}
{\sc J. Wolf}, {\em On the local
  regularity of suitable weak solutions to the generalized navier-stokes
  equations}, Ann Univ Ferrara, 61 (2015), pp.~149--171.

\bibitem{Wolf2016}
\leavevmode\vrule height 2pt depth -1.6pt width 23pt, {\em On the local
  pressure of the navier-stokes equations and related systems}, to appear in
  Adv. Diff. Equs.,  (2016).

\end{thebibliography}

\end{document}